\documentclass[11pt,a4paper]{article}

\usepackage[utf8]{inputenc}
\usepackage[T1]{fontenc}
\usepackage{lmodern}
\usepackage[margin=2.6cm]{geometry}
\usepackage{amsmath,amssymb,amsthm,mathtools}
\usepackage{booktabs}
\usepackage{graphicx}
\graphicspath{{figures/}}
\usepackage{enumitem}
\usepackage[numbers,sort&compress]{natbib}
\usepackage[colorlinks=true,linkcolor=blue!50!black,citecolor=green!40!black,urlcolor=blue!50!black]{hyperref}
\usepackage{xcolor}

\theoremstyle{plain}
\newtheorem{theorem}{Theorem}[section]
\newtheorem{lemma}[theorem]{Lemma}
\newtheorem{proposition}[theorem]{Proposition}
\newtheorem{corollary}[theorem]{Corollary}
\theoremstyle{definition}

\newtheorem{example}[theorem]{Example}
\theoremstyle{remark}
\newtheorem{remark}[theorem]{Remark}

\newcommand{\R}{\mathbb{R}}
\newcommand{\eps}{\varepsilon}
\newcommand{\ce}[1]{\mathrm{#1}}
\newcommand{\ip}[2]{\langle #1,#2\rangle}
\DeclareMathOperator{\rank}{rank}
\DeclareMathOperator{\sign}{sign}
\DeclareMathOperator{\relint}{relint}

\title{Canard explosions in mass-action networks\\ with one slow direction}
\author{Miruna-Mihaela Beldiman\thanks{Institute of Mathematical Statistics
and Applied Mathematics of the Romanian Academy, Bucharest, Romania.
ORCID 0000-0001-9947-7594.}
\and
Andrei-Dan Halanay\thanks{University of Bucharest, Bucharest, Romania.
ORCID 0009-0005-5481-8021. Corresponding author:
\href{mailto:andrei.halanay@unibuc.ro}{andrei.halanay@unibuc.ro}.}}
\date{}

\begin{document}
\maketitle

\begin{abstract}
We study canard explosions in mass-action networks with one slow
direction. For any number of species and any rank of the fast subnetwork,
the derivatives that enter the Krupa--Szmolyan conditions and coefficient
are rational functions of the fluxes and the reciprocal concentrations.
For a fast subnetwork of rank one,
the dynamics is planar on stoichiometric classes and these forms are linear
in the fluxes, extending our planar flux criterion; the criticality is an
explicit cubic. When one species is fast, these forms decompose over the
slow species, weighted by their reciprocal concentrations. The
elementary-mode criterion then carries over, and the
criticality depends on the composition of the slow pool, so one network
can explode in either way. The obstruction that the slow
flux cannot vanish on the fold is characterised by a Farkas certificate and
a yield diagram decided on configurations of at most three reactions.
Under four conditions on the critical curve, a canard point with $A\neq0$
unfolded by a rate constant gives the full canard explosion; with one fast
species and a fast equation affine in the slow variable, these follow
from sign conditions on the reactions and injectivity on the class. For fast subnetworks of higher rank, a centre-manifold reduction adds
corrections through the group inverse of the fast Jacobian; in the
Oregonator they make both canard points change criticality between its
planar reduction and the full model. All results are proved; exact computations and
simulations serve as independent checks.
\end{abstract}



\section{Introduction}\label{sec:intro}

Oscillations in chemical and biochemical reaction systems often start
abruptly: as a parameter crosses a threshold, small oscillations around a
steady state turn into large relaxation oscillations within a parameter
interval so small that the transition looks like a bifurcation of its own
\citep{PengGasparShowalter1991}. This canard explosion has been found in
models of the Belousov--Zhabotinsky reaction
\citep{BronsBarEli1991,RotsteinKopellZhabotinskyEpstein2003}, and slow--fast
analyses of the same kind underlie models of the peroxidase--oxidase
reaction \citep{Olsen1983,KuehnSzmolyan2015}, of intracellular calcium
dynamics \citep{HarveyKirkWechselbergerSneyd2011} and of the embryonic cell
cycle \citep{KosiukSzmolyan2016}. The kind of onset matters. In a
supercritical explosion the amplitude grows through stable small cycles
before the jump; in a subcritical one the small cycles are unstable, a
stable steady state coexists with the large oscillation, and the system
shows hysteresis. A reaction scheme is usually known before its rate
constants, which are often uncertain. We therefore ask a structural
question: which reaction networks can show a canard explosion, for some
rate constants and separation of time scales, and of which kind, answered
from the network and its steady-state fluxes rather than from a particular
choice of parameters.

Mathematically, a canard explosion is the transition, within an
exponentially small parameter interval, from small limit cycles born in a
Hopf bifurcation to relaxation oscillations. It occurs in slow--fast systems at a canard point,
a fold of the critical manifold at which the slow flow passes from the
attracting to the repelling branch. The phenomenon was discovered with
nonstandard analysis \citep{BenoitCallotDiener1981} and explained
geometrically by \citet{KrupaSzmolyan2001fold,KrupaSzmolyan2001relax},
whose coefficient $A$ decides whether the explosion is supercritical or
subcritical; see \citet{Kuehn2015}, \citet{Wechselberger2020} and
\citet{DeMaesschalckDumortierRoussarie2021} for the theory and
\citet{DeMaesschalckDoanWynen2021} for an intrinsic criterion for the
criticality. Chemical reactions were among the first sources of canards
\citep{BronsBarEli1991,BronsSturis2001,GucwaSzmolyan2009,KuehnSzmolyan2015,HasanKrauskopfOsinga2017},
but each analysis concerns one model, in its own coordinates.

In \citet{BeldimanHalanay2026} we answered the structural question
for networks of two species, one fast and one slow. In flux coordinates the
answer is a criterion of three linear equations, an inequation and a sign
condition, the criticality is the sign of a cubic form, and a network has a
canard point iff its elementary flux modes contain a sign-meeting pair and
a mode carrying a fold. The present paper extends this theory to networks
with any number of species and one slow direction. Most of the paper
concerns fast subnetworks of rank one. Conservation laws then make the
dynamics planar on each stoichiometric class, but the planar coordinates mix
the species, and three new features appear: the fast direction may combine several species,
the slow direction may involve several slow species, and we take up the
global dynamics, which the planar paper left aside.
Slow--fast splittings of this kind arise from scalings of rate constants
and concentrations
\citep{GoekeWalcherZerz2015,FeliuKruffWalcher2020,FeliuWalcherWiuf2022,Kruff2021},
and the structural setting is that of chemical reaction network theory
\citep{Feinberg2019}: classifications of small networks
\citep{JoshiShiu2017,BanajiBoros2023}, conditions on Jacobian determinants
\citep{ConradiFeliuMinchevaWiuf2017}, inheritance of bifurcations
\citep{Banaji2018,BanajiBorosHofbauer2025} and injectivity
\citep{CraciunFeinberg2005,ShinarFeinberg2012}. Our coordinates, fluxes and
reciprocal concentrations, are the convex coordinates of stoichiometric
network analysis \citep{Clarke1980}, which have been used to locate Hopf
bifurcations of reaction networks
\citep{GatermannEiswirthSensse2005,ErramiEtAl2015}; here they serve to
separate the geometry of the critical curve, which they describe exactly,
from the choice of rate constants.

Our first result is the flux criterion for this class
(Theorem~\ref{thm:flux}). At a point $x$ with fluxes $v$, all derivatives
of the planar system that enter the conditions of Krupa and Szmolyan and
their coefficient $A$ are linear forms in $v$ whose coefficients are
polynomials in $h=1/x$, built from the relative changes of
the concentrations along the fast and the slow directions. Each of these
forms has a meaning in the phase plane of the class (Table~\ref{tab:dict}
and Figure~\ref{fig:geometry}), and the cubic that gives the criticality
compares the curvature of the slow flux along the critical curve with the
skewness of the fold. Canard points are
characterised as in the planar case, with forms that now depend on the
point. When a single species $\ce X$ is fast, the dependence collapses
(Section~\ref{sec:one}): the canard cone is independent of the point, the
sensitivity to the slow variable and the criticality cubic are sums over
the slow species $i$ of the planar quantities of the pairs $(\ce X,i)$
weighted by $U_i/x_i$, and the mode criterion carries over. A consequence
without planar analogue is that the criticality depends on the composition
of the slow pool: even a network carried by one elementary mode can have a
supercritical singular Hopf bifurcation at some concentrations and a
subcritical one at others.

The second result concerns sign-locking, the obstruction that the slow
flux has a fixed sign at every fold point for all rate constants
(Section~\ref{sec:lock}); among the randomly sampled networks that fold, it
is the obstruction we met most often. We show that this happens exactly
when the slow flux minus a combination of the fold equations has
nonnegative coefficients, and that this is decided by a yield diagram in
which every reaction is drawn at its order and its slow yield, by
configurations of at most three reactions. With the fold criterion and the
mode criterion this yields a complete hierarchy of obstructions. The
certificate has a dynamical meaning: in a sign-locked network the slow flow
crosses every fold in the same direction, so for no rate constants is an
oscillation born in a singular Hopf bifurcation near a fold, once $\eps$ is
small (Proposition~\ref{prop:lockdyn}).

The third result is global (Section~\ref{sec:global}). Under four
conditions on the critical curve, which we reduce to the hypotheses of
\citet{KrupaSzmolyan2001relax}, a canard point with $A\neq0$ that is
unfolded by a parameter gives the full canard explosion, and the distance
between the Hopf point and the explosion is again given by the fluxes
(Corollary~\ref{cor:location}). For networks with one fast species whose fast
equation is affine in the slow variable these conditions are a finite
computation at given rate constants (Proposition~\ref{prop:decide}). It
applies, for instance, to the two-variable Oregonator of
\citet{TysonFife1980}, which for every value $q<0.0797$ of its parameter $q$
undergoes the full explosion at both of its canard points, with a change of
criticality at one of them as $q$ varies (Proposition~\ref{prop:oreg2}). For networks
with one fast species whose fast equation is affine in the slow variable,
four sign conditions on the reactions give three of the four conditions
whenever the second fold lies in the class, and two more, one of which is
injectivity on the class, make every canard point of the network with
$\mathcal Q_h\neq0$ explode, whatever the rate constants, when $\rank B=3$.
A six-reaction network on three species illustrates both explosion
scenarios in one network, the criticality being set by the ratio of the
slow species. Finally, for a fast direction that mixes two species, the
fold condition is bilinear in fluxes and reciprocal concentrations, and
fold equilibria are confined to a window of concentration ratios computed
from elementary modes (Section~\ref{sec:mixing}).

The fourth result removes the assumption of rank one
(Section~\ref{sec:rank}). For a fast subnetwork of any rank $s$, a reduction
to a centre manifold of the fold gives a planar slow--fast system whose
Krupa--Szmolyan coefficients are again forms in the fluxes, homogeneous of
degree one, now with corrections through the group inverse of the fast
Jacobian, which describe how the transverse fast species are slaved to the
fold direction, that is, determined by it on the centre manifold
(Theorem~\ref{thm:rank}). The fold condition becomes in
general a form of degree $s$ in the fluxes, so the polyhedral theory of the earlier sections,
with its elementary modes and certificates, is specific to rank one, while
the flux criterion and the formula for the criticality hold for every rank.
The three-variable Oregonator has a fast subnetwork of rank two. Both of its
canard points change criticality between its planar reduction and the
three-variable model, and Theorem~\ref{thm:rank} shows that the time-scale
ratio of the two fast species acts on the criticality only through the
corrections. The first Lyapunov coefficient confirms the signs, and
simulation the subcritical explosion at one of the points
(Example~\ref{ex:oregonator}).

Sections~\ref{sec:flux} and~\ref{sec:one} extend the criteria of
\citet{BeldimanHalanay2026} to the rank-one class. The paper is
self-contained: the two results of \citet{BeldimanHalanay2026} that the
proofs use, the formula for the coefficient $A$ and the mode criterion, are
proved here (Appendix~\ref{app:A} and the proof of Theorem~\ref{thm:modes}).
What is new is the dependence of the criticality on
the composition of the slow pool (Proposition~\ref{prop:crit} and
Example~\ref{ex:onemode}; a complete enumeration of a three-species
universe finds $560$ of its $1\,302$ single-mode canard networks with both
criticalities), the theory
of sign-locking (Section~\ref{sec:lock}), which is new also for two
species, the global explosion (Section~\ref{sec:global}), the fold
window of a mixing fast direction (Section~\ref{sec:mixing}) and the
criterion for fast subnetworks of higher rank, with its application to the
Oregonator (Section~\ref{sec:rank}). All results
concern the limit $\eps\to0$; Remark~\ref{rem:BT} estimates how small
$\eps$ must be. Every theorem, proposition and corollary is proved. The
computations, by exact linear programming, exact real algebra, first
Lyapunov coefficients and simulation, are independent checks of the proofs
and of the examples, reported in the supplementary materials; frequencies in
samples of networks, as in Section~\ref{sec:lock}, describe the samples and
are used in no proof.

Section~\ref{sec:setting} sets up the class of networks and
Section~\ref{sec:flux} proves the flux criterion. Section~\ref{sec:one}
treats one fast species, Section~\ref{sec:lock} sign-locking,
Section~\ref{sec:global} the global explosion and Section~\ref{sec:mixing}
a mixing fast direction, and Section~\ref{sec:rank} fast subnetworks of
higher rank. Section~\ref{sec:disc} discusses open problems, and
Appendix~\ref{app:A} derives the coefficient $A$ of a planar canard point.

\section{Setting}\label{sec:setting}

Species $X_1,\dots,X_n$ have concentrations $x\in\R^n_{>0}$, and reaction
$r$ has reactant complex $y_r\in\mathbb Z^n_{\ge0}$ and mass-action rate
$v_r=k_rx^{y_r}$. After a scaling, the change of each reaction splits into
a fast and a slow part, and we assume that all fast parts are parallel to
one vector $w$:
\begin{equation}\label{eq:sys}
  \dot x=w\sum_rL_rv_r(x)+\eps\sum_rN^1_rv_r(x)
\end{equation}
in fast time, with $L_r\in\R$ and $N^1_r\in\R^n$; $R$ is the number of
reactions. We assume:
\begin{enumerate}[label=(H\arabic*),nosep]
  \item \emph{rank one:} the fast stoichiometric subspace is $\R w$;
  \item \emph{one slow direction:} $\mathcal S=\R w+\operatorname{span}\{N^1_r\}$ has
        dimension two.
\end{enumerate}
For $n>2$, (H2) means that the network has at least $n-2$ conservation
laws, and exactly $n-2$ when the reaction vectors span $\mathcal S$.
Choose $U\in\mathcal S$ and a linear form $\varpi$ with $\varpi w=0$ and
$\varpi U=1$, and let $g_r=\varpi N^1_r$. On the class through $x^*$, the
coordinates $x=x^*+zw+uU$ turn \eqref{eq:sys} into the planar slow--fast
system
\begin{equation}\label{eq:planar}
  z'=F(z,u)+\eps\chi(z,u),\qquad u'=\eps G(z,u),\qquad
  F=\sum_rL_rv_r,\quad G=\beta=\sum_rg_rv_r,
\end{equation}
where $\chi$ collects the slow parts along $w$; it affects neither the
conditions below nor the coefficient $A$. Three cases guide the paper.
For $n=2$, $w=e_{\ce X}$ and $U=e_{\ce Y}$, this is the system of
\citet{BeldimanHalanay2026}. If $w=e_{\ce X}$ and $n>2$, a single species is
fast and the slow direction may involve several conserved slow species, as
in $\ce Y\rightleftharpoons\ce Z$ with $\ce Y+\ce Z$ conserved. If $w$ has
several nonzero entries, as for fast conversions $\ce X\rightleftharpoons\ce Z$,
the fast direction mixes species.

The splitting is a property of a regime of rate constants and
concentrations, not of the network. A network can have several regimes of
the form~\eqref{eq:sys}, or none, and our results concern the canard points
within one of them. Such regimes are found by the scaling methods cited in
the introduction; in Example~\ref{ex:onemode} the slow parts come from
abundant slow species, and in Example~\ref{ex:both} from slow rate
constants. If the fast parts of the reactions are not parallel to one
vector, (H1) fails and the fast subnetwork has rank two or more
(Section~\ref{sec:rank}).

The coefficient $A$ of Theorem~\ref{thm:flux}(b) below does not depend on
the choices made here: it is unchanged when $w$ or $U$ is multiplied by a
positive number, when $U$ is replaced by $U+cw$, and under the reflections
$z\mapsto-z$, $u\mapsto-u$. Conditions away from the canard point can depend
on $U$: replacing $U$ by $U+cw$ replaces $F_u$ by $F_u+cF_z$, which agrees
with $F_u$ only on the fold. We state such conditions, for instance (G1) in
Section~\ref{sec:global}, for a fixed choice of $U$; it suffices that they
hold for one choice.

\section{The flux criterion}\label{sec:flux}

Mass-action rates are monomials, and the derivative of a monomial along a
direction $a$ is $\partial_a x^{y}=x^{y}\ip{y}{a/x}$: it is the rate times
a linear form in the relative changes $a_i/x_i$. So every derivative of the
vector field at a point is a combination of the fluxes $v_r$ with
coefficients that are polynomials in the reciprocal concentrations, and,
since every flux vector $v>0$ is realised at every point by suitable rate
constants, fluxes and points can be prescribed independently. This is why we
use the coordinates $(h,v)$, with $h=1/x$ componentwise: they are the convex
coordinates of \citet{Clarke1980}, and in them the local geometry of the
critical curve at a point (fold, curvature, position of the slow nullcline)
is linear in $v$. Define the log-directions $\delta=w\circ h$ and
$\kappa=U\circ h$, the relative changes of the concentrations along the fast
and the slow direction.

\begin{lemma}\label{lem:deriv}
Along $x+zw+uU$, at $z=u=0$, with powers of vectors taken componentwise,
\begin{align*}
  \partial_zv_r&=v_r\ip{y_r}{\delta}, &
  \partial_uv_r&=v_r\ip{y_r}{\kappa},\\
  \partial_z^2v_r&=v_r\big(\ip{y_r}{\delta}^2-\ip{y_r}{\delta^2}\big), &
  \partial_z\partial_uv_r&=v_r\big(\ip{y_r}{\delta}\ip{y_r}{\kappa}-\ip{y_r}{\delta\kappa}\big),
\end{align*}
and
$\partial_z^3v_r=v_r\big(\ip{y_r}{\delta}^3-3\ip{y_r}{\delta}\ip{y_r}{\delta^2}+2\ip{y_r}{\delta^3}\big)$.
\end{lemma}

\begin{proof}
$\partial_zv_r=v_rT_1$ with $T_k=\sum_iy_{ri}w_i^k/(x_i+zw_i+uU_i)^k$, and
$\partial_zT_k=-kT_{k+1}$; differentiate, and mix in $u$ in the same way.
\end{proof}

Accordingly we define, for $h>0$ and $v\in\R^R$, the weights
$c^d_r=\ip{y_r}{\delta}$, $c^e_r=\ip{y_r}{\kappa}$,
$c^{dd}_r=\ip{y_r}{\delta}^2-\ip{y_r}{\delta^2}$,
$c^{de}_r=\ip{y_r}{\delta}\ip{y_r}{\kappa}-\ip{y_r}{\delta\kappa}$ and
$c^{ddd}_r$ as in the lemma, and the linear forms
\begin{gather*}
  F_0=\textstyle\sum_rL_rv_r,\quad F_z=\sum_rL_rc^d_rv_r,\quad
  F_{zz}=\sum_rL_rc^{dd}_rv_r,\quad F_{zzz}=\sum_rL_rc^{ddd}_rv_r,\\
  F_u=\textstyle\sum_rL_rc^e_rv_r,\quad F_{zu}=\sum_rL_rc^{de}_rv_r,\quad
  G_0=\sum_rg_rv_r,\quad G_z=\sum_rg_rc^d_rv_r,\\
  G_u=\textstyle\sum_rg_rc^e_rv_r,\quad G_{zz}=\sum_rg_rc^{dd}_rv_r .
\end{gather*}
By Lemma~\ref{lem:deriv}, if $v$ are the fluxes at $x=1/h$, these are the
corresponding derivatives of $F$ and $G$ in~\eqref{eq:planar}. Let
\[
  B_h=\begin{pmatrix}L\\ L\circ c^d(h)\\ g\end{pmatrix},\qquad
  K_h=\{v>0:B_hv=0\},
\]
the \emph{canard cone} at $h$. Here $B_h$ is a $3\times R$ matrix. When
$K_h$ is nonempty it is a relatively open polyhedral cone of dimension
$R-\rank B_h$, and by Stiemke's theorem $K_h\neq\emptyset$ iff there is no
$\lambda\in\R^3$ with $\lambda^{\mathsf T}B_h\ge0$ and
$\lambda^{\mathsf T}B_h\neq0$, which is a linear program. Let
\begin{equation}\label{eq:Qh}
  \mathcal Q_h(v)=(G_zF_{zu}-G_uF_{zz})F_{zz}+F_u(G_{zz}F_{zz}-G_zF_{zzz}).
\end{equation}
A \emph{canard point} of \eqref{eq:planar} is a point with $F=F_z=G=0$,
$F_{zz}\neq0$ and $F_uG_z<0$. The inequations contain the nondegeneracy
conditions~(3.2) of \citet{KrupaSzmolyan2001relax} other than
$g_\lambda\neq0$, and the sign condition is the orientation of their
canonical form~(3.3) and of their assumption (A4$'$), recalled in
Section~\ref{sec:global}: the slow flow crosses
the fold from the attracting to the repelling branch. A parameter
\emph{unfolds} the canard point if the derivative of $G$ along the
continuation of the fold point does not vanish, which is the remaining
condition $g_\lambda\neq0$.

\begin{theorem}[Flux criterion]\label{thm:flux}
Assume (H1)--(H2).
\begin{enumerate}[label=(\alph*),nosep]
  \item For given rate constants, $x$ is a canard point iff its fluxes
        satisfy $v\in K_h$, $F_{zz}(v,h)\neq0$ and $F_u(v,h)G_z(v,h)<0$, with
        $h=1/x$. Hence some rate constants give a canard point iff these
        conditions have a solution $(h,v)\in\R^n_{>0}\times\R^R_{>0}$.
  \item At a canard point
        $A=2\mathcal Q_h(v)/\big(F_{zz}^2\sqrt{-F_uG_z}\big)$, so
        $\sign A=\sign\mathcal Q_h(v)$.
  \item Fix the stoichiometric class. If $\rank B_h=3$, the rate constants
        for which the class has a canard point near $x$ form a smooth
        hypersurface near $k$, and some single rate constant unfolds the
        canard point. If $g$ lies in the span of $L$ and $L\circ c^d(h')$
        for all $h'$ near $h$, which with one fast species is the case
        $\rank B_h<3$ (at a canard point the rows $L$ and $L\circ c^d$ are
        independent, since $L\circ a=cL$ would give $n_{20}=c(c-1)Lv=0$), the set of such rate constants has nonempty
        interior.
\end{enumerate}
\end{theorem}

\begin{proof}
(a) Every $v>0$ is realised at every $x$ by $k_r=v_r/x^{y_r}$, and
Lemma~\ref{lem:deriv} identifies the conditions. (b) The system on a class
is planar, of the form (3.1) of \citet{KrupaSzmolyan2001relax}, and
Lemma~\ref{lem:A} gives $A$; the term $\eps\chi$ enters only the coefficient
$a_1$, which does not appear in $A$. (c) On the class, the Jacobian of
$(F,\partial_zF)$ in $(z,u)$ at the canard point has determinant
$-F_uF_{zz}\neq0$, so the fold point continues smoothly in $k$ as $p_F(k)$,
and near $k$ the rate constants with a canard point on the class form the
zero set of $\gamma(k)=G(p_F(k),k)$. The derivative of $(F,\partial_zF,G)$
with respect to $k$ at fixed $x$ is $B_h\operatorname{diag}(x^{y_r})$. If
$\rank B_h=3$, there is $\delta k$ with $D_k(F,\partial_zF)\,\delta k=0$ and
$D_kG\,\delta k=1$; along $\delta k$ the fold point does not move to first
order, so $\nabla_k\gamma\cdot\delta k=1$. Hence $\nabla_k\gamma\neq0$, the
zero set is a smooth hypersurface, and a rate constant $k_r$ with
$\partial\gamma/\partial k_r\neq0$ unfolds the canard point. In
the second case $G$ is a combination of $F$ and $F_z$ at every point near
$x$, so $G=0$ at the continued fold point and the remaining conditions are
open.
\end{proof}

\begin{figure}[ht]
\centering
\includegraphics[width=\textwidth]{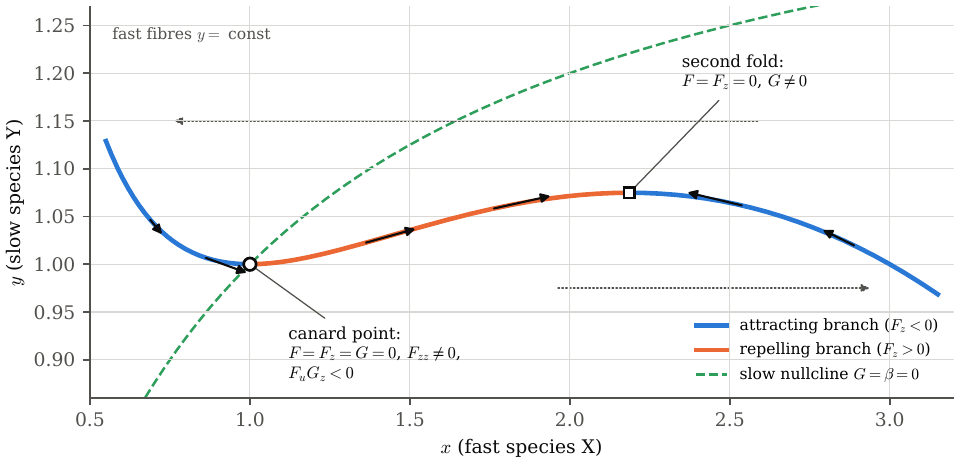}
\caption{The geometry of a canard point, for the network of
Example~\ref{ex:both} with $s=5$, $m=1$ and $y/z=2$ at the canard point,
in the coordinates $(x,y)$ of the class. The critical curve $F=0$ has a
minimum, the canard point, where the slow nullcline crosses it, and a
maximum, a jump point. Arrows on the curve show the slow flow, dotted arrows
the fast fibres.}
\label{fig:geometry}
\end{figure}

\paragraph{Geometric meaning.}
Each condition of Theorem~\ref{thm:flux} is a statement about the phase plane
of the class (Table~\ref{tab:dict} and Figure~\ref{fig:geometry}): $F=0$
places the point on the critical curve $C_0$, $F_z=0$ makes it a fold, and the
sign of $F_z$ distinguishes attracting from repelling branches; as a graph
$u=\varphi(z)$, the critical curve has $\varphi''=-F_{zz}/F_u$ at the fold.
The condition $G=0$ puts the slow nullcline through the fold, and $F_uG_z<0$
says that the slow flow carries trajectories along the attracting branch
into the fold and on along the repelling branch; for $F_uG_z>0$ the
equilibrium at the fold is a saddle. The canard cone $K_h$ is thus the set of
flux vectors at $x$ for which the slow nullcline meets the critical curve at
a fold. The criticality has an intrinsic meaning as well. Parametrise $C_0$
near the fold by $s=F_z$, let $\mathcal G(s)$ be the slow flux multiplied by
the transversality factor $F_uF_{zz}-F_zF_{zu}$ along $C_0$, and let
$\mathcal A=F_{zzz}/F_{zz}^2$ be the skewness of the fold.
\citet{DeMaesschalckDoanWynen2021} give an intrinsic criterion; in the chart
$(z,u)$, with $Z=\partial_z$ and the Euclidean metric, it says that the
slow--fast Hopf point is supercritical if
$\sigma=\tfrac12\mathcal G''(0)-\mathcal G'(0)\mathcal A<0$ and subcritical
if $\sigma>0$. At a canard point $\mathcal G'(0)=F_uG_z$ and
$\sigma=\mathcal Q_h/(2F_{zz}^2)$ (Section~\ref{sm:rank}). So the cubic
$\mathcal Q_h$ compares the curvature of the slow-flux profile along the
critical curve with the skewness of the fold, both written as forms in the
fluxes.

\begin{table}[ht]
\centering
\caption{The objects of the flux criterion and their meaning in the phase
plane of a class, for a fast subnetwork of rank one.}\label{tab:dict}
\small
\begin{tabular}{@{}p{0.34\textwidth}p{0.6\textwidth}@{}}
\toprule
object & meaning\\
\midrule
$F=\sum_rL_rv_r=0$ & the point lies on the critical curve $C_0$\\
$F_z=0$; sign of $F_z$ & fold of $C_0$; attracting ($F_z<0$) or repelling ($F_z>0$) branch\\
$F_{zz}\neq0$, $F_u$ & quadratic fold; $\varphi''=-F_{zz}/F_u$ for $C_0$ as a graph $u=\varphi(z)$\\
$G=\beta=0$ & the slow nullcline passes through the point\\
$F_uG_z<0$ & the slow flow passes through the fold from the attracting to the repelling branch\\
canard cone $K_h$ & flux vectors at $x$ for which the slow nullcline meets $C_0$ at a fold\\
cubic $\mathcal Q_h$ & curvature of the slow flux along $C_0$ against the skewness of the fold\\
fold-flux cone $\Phi$ (Section~\ref{sec:lock}) & flux vectors of all fold points, for all rate constants\\
separating line of the yield diagram & the slow flow crosses every fold in the same direction, for all rate constants (Proposition~\ref{prop:lockdyn})\\
elementary modes & minimal sets of reactions that balance the fast and slow fluxes at a fold\\
\bottomrule
\end{tabular}
\end{table}

The hypothesis on nearby $h'$ cannot be weakened to $\rank B_h<3$ when the
fast direction mixes species. For the fast reactions $\ce X\to\ce Z$,
$2\ce X+\ce Z\to3\ce X$ and the slow reactions $\ce X\to\ce Y$,
$\ce Y\to\ce Z$, with rate constants $k_1,\dots,k_4$ in this order,
$w=e_{\ce X}-e_{\ce Z}$ and $\ce X+\ce Y+\ce Z$
conserved, $F_z$ vanishes at fluxes with $Lv=0$ exactly when $x=z$, where
$\rank B_h=2$; the point $(1,1,1)$ with unit fluxes is a canard point. In
general, near a canard point the rate constants with a canard point on the
class form the zero set of $\gamma(k)=G(p_F(k),k)$, where $p_F(k)$ is the
continued fold point, and
\[
  \nabla_k\gamma=D_kG-(G_z,G_u)\,J^{-1}D_k(F,F_z),
\]
with $J$ the Jacobian of $(F,F_z)$ in $(z,u)$; the set is a smooth
hypersurface wherever $\nabla_k\gamma\neq0$. The case $\rank B_h=3$ of
Theorem~\ref{thm:flux}(c) is the one in which a perturbation of $k$ can
leave $p_F$ fixed to first order. In the example, the fold lies at
$x=z=\sqrt{k_1/k_2}$, and with $y=T-2x$ on the class of total $T$ the
condition $G=0$ reads $k_3x=k_4(T-2x)$; so $\partial\gamma/\partial k_3=x\neq0$,
and on a fixed class the rate constants with a canard point form a
hypersurface.

For $n=2$, $w=e_{\ce X}$ and $U=e_{\ce Y}$, the forms are the moment forms of
\citet{BeldimanHalanay2026} divided by powers of $x$ and $y$:
$F_z=n_{10}/x$, $F_{zz}=n_{20}/x^2$, $F_u=n_{01}/y$, $G_z=m_{10}/x$ and so
on, and $\mathcal Q_h=Q/(x^4y)$. For a mixing fast direction we use two networks on
$\ce X,\ce Y,\ce Z$ with $\ce X+\ce Y+\ce Z$ conserved, $w=e_{\ce X}-e_{\ce Z}$
and $U=e_{\ce Y}-e_{\ce Z}$. Network~A has the fast reactions
$\ce Z\to\ce X$, $\ce X\to\ce Z$, $\ce Z+2\ce X\to3\ce X$,
$\ce X+\ce Y\to\ce Z+\ce Y$ and the slow reactions $\ce X\to\ce Y$,
$\ce Y\to\ce Z$; network~B has the reactions $\ce Z\to\ce Y$,
$\ce Z\to\ce X$, $\ce X\to\ce Z$, $\ce X+\ce Y\to\ce X+\ce Z$,
$\ce X+\ce Y+\ce Z\to2\ce X+\ce Y$, $\ce Z+2\ce X\to3\ce X$, of which the
first and the fourth are slow. They have canard points at $72$ and $20$ of
the $216$ points of the grid $\{1,\dots,6\}^3$, respectively, and admit both
criticalities. As independent checks, the identities, Lemma~\ref{lem:deriv}
and Theorem~\ref{thm:flux}(b) were also verified by exact computation and
against the first Lyapunov coefficient of the Hopf bifurcation on the class
(Section~\ref{sm:flux} of the supplementary materials).

\section{One fast species}\label{sec:one}

Let $w=e_{\ce X}$. Replacing $U$ by $U-U_{\ce X}e_{\ce X}$ we may assume
$U_{\ce X}=0$; the \emph{slow species} are those with $U_i\neq0$. Write
$a_r=y_{r\ce X}$ and, for each species $i$, define the moment forms of the
pair $(\ce X,i)$ as in \citet{BeldimanHalanay2026}, with the falling
factorial $(a)_k=a(a-1)\cdots(a-k+1)$,
\begin{gather*}
  n_{k0}=\sum_rL_r(a_r)_kv_r,\qquad m_{k0}=\sum_rg_r(a_r)_kv_r,\\
  n^{(i)}_{01}=\sum_rL_ry_{ri}v_r,\qquad n^{(i)}_{11}=\sum_rL_ra_ry_{ri}v_r,\qquad
  m^{(i)}_{01}=\sum_rg_ry_{ri}v_r,
\end{gather*}
and the planar cubic
$Q^{(i)}=(m_{10}n^{(i)}_{11}-m^{(i)}_{01}n_{20})n_{20}+n^{(i)}_{01}(m_{20}n_{20}-m_{10}n_{30})$.

\begin{proposition}[Slow-pool decomposition]\label{prop:sup}
$K_h=K=\{v>0:\ \sum_rL_rv_r=\sum_rL_ra_rv_r=\sum_rg_rv_r=0\}$ does not depend
on the point, $F_{zz}=h_{\ce X}^2n_{20}$ and $G_z=h_{\ce X}m_{10}$, and
\[
  F_u=\sum_iU_ih_in^{(i)}_{01},\qquad
  \mathcal Q_h=h_{\ce X}^4\sum_iU_ih_i\,Q^{(i)} .
\]
\end{proposition}

\begin{proof}
$\delta=h_{\ce X}e_{\ce X}$ and $\delta\circ\kappa=0$, so $c^d_r=h_{\ce X}a_r$,
$c^{dd}_r=h_{\ce X}^2(a_r)_2$, $c^{ddd}_r=h_{\ce X}^3(a_r)_3$,
$c^e_r=\sum_iU_ih_iy_{ri}$ and $c^{de}_r=h_{\ce X}a_rc^e_r$. Substituting into
\eqref{eq:Qh}, every term is $h_{\ce X}^4$ times a sum over $i$ of $U_ih_i$
times the corresponding term of $Q^{(i)}$.
\end{proof}

The decomposition itself is elementary: $F_u$ and $\mathcal Q_h$ are linear
in the weights $U_ih_i$ of the slow species. Its consequence is not: the
criticality, and not only the location, of a canard point depends on the
composition of the slow pool (Proposition~\ref{prop:crit}). The forms $n_{k0}$ and $m_{k0}$,
among them the fold form $n_{20}$ and the sensitivity $m_{10}$, are shared,
and only $n^{(i)}_{01}$, $n^{(i)}_{11}$ and $m^{(i)}_{01}$ depend on $i$. An
\emph{elementary mode} of $B=[L;\,L\circ a;\,g]$ is a nonnegative vector in
its kernel of minimal support, that is, an extreme ray of the closure of
$K$. For an elementary mode $e$ and a slow species $i$, put
$p_i(e)=(U_in^{(i)}_{01}(e),m_{10}(e))$; two modes, possibly equal, are
\emph{sign-meeting for $i$} if the cone spanned by $p_i(e)$ and $p_i(e')$ meets
$\Omega=\{\varphi\psi<0\}$.

\begin{theorem}[Mode criterion]\label{thm:modes}
A network with $K\neq\emptyset$ and one fast species has a canard point iff
its elementary modes include a mode with $n_{20}\neq0$ and, for some slow
species $i$, a possibly equal pair that is sign-meeting for $i$. Every such network
contains a subnetwork with a canard point that is carried by at most three modes,
hence has at most twelve reactions, independently of the number of species
and of the molecularity.
\end{theorem}

\begin{proof}
For fixed $v$ the sign condition reads $m_{10}\sum_iU_ih_in^{(i)}_{01}<0$.
Since $h$ ranges over $\R^n_{>0}$, it holds for some $h$ iff
$p_i(v)\in\Omega$ for some $i$: a negative sum has a negative term, and if
$U_in^{(i)}_{01}m_{10}<0$, the sum is negative when $h_j$ is small for the
other slow species $j$. Since also $F_{zz}=h_{\ce X}^2n_{20}$, by
Theorem~\ref{thm:flux}(a) the network has a canard point iff some $v\in K$
and some $i$ satisfy $p_i(v)\in\Omega$ and $n_{20}(v)\neq0$.

Suppose $v$ is such. The closure $\{v\ge0:Bv=0\}$ of $K$ is a pointed
polyhedral cone with relative interior $K$, so $v=\sum_jc_je_j$ over all
elementary modes with all $c_j>0$: $v-t\sum_je_j$ lies in the cone for small
$t>0$, and is a nonnegative combination of the modes. Then $n_{20}(e_j)\neq0$ for some $j$,
and $p_i(v)=\sum_jc_jp_i(e_j)$ lies in $\Omega$ and, by Carathéodory's
theorem for cones in $\R^2$, in the cone spanned by at most two of the
$p_i(e_j)$, which are therefore sign-meeting for $i$. Conversely, let $e,e'$
be sign-meeting for $i$, with $p_i(se+s'e')\in\Omega$ for some $s,s'\ge0$,
and let $n_{20}(z)\neq0$ for a mode $z$. Since $\Omega$ is open,
$w=se+s'e'+\gamma z$ has $p_i(w)\in\Omega$ for small $\gamma>0$, and
$n_{20}(w)\neq0$ for all but at most one $\gamma$. Then $v=w+\delta u$, with
$u\in K$ and $\delta>0$ small, lies in $K$ and satisfies both conditions.

For the second statement, let $W$ be the union of the supports of $e$,
$e'$ and $z$. These are elementary modes of the subnetwork $W$, and
$e+e'+z$ is a positive flux of $W$, so the first statement applies to $W$.
A minimal support $T$ of a kernel vector of $B$ has $\rank B_T=|T|-1\le3$,
so every mode has at most four reactions and $W$ at most twelve.
\end{proof}

\begin{proposition}\label{prop:crit}
At a canard point, $\sign A=\sign\sum_iU_ih_iQ^{(i)}(v)$, a linear form in
$h$. With two or more slow species both signs can occur for the same
fluxes, even when $K$ is a single ray.
\end{proposition}

\begin{proof}
The first statement follows from Theorem~\ref{thm:flux}(b) and
Proposition~\ref{prop:sup}, since $h_{\ce X}^4>0$. Example~\ref{ex:onemode}
proves the second.
\end{proof}

\begin{example}[One mode, both criticalities]\label{ex:onemode}
Take $0\to\ce X$, $\ce Z\to\ce Y$, $\ce X+\ce Y+\ce Z\to2\ce Z$ and
$2\ce X\to3\ce X$, with $\ce X$ fast, $\ce Y+\ce Z$ conserved and
$U=e_{\ce Y}-e_{\ce Z}$. Its only mode is $e=(1,2,2,1)$, with $n_{20}=2$,
$m_{10}=-2$, $n^{(\ce Y)}_{01}=n^{(\ce Z)}_{01}=-2$, $Q^{(\ce Y)}=16$ and
$Q^{(\ce Z)}=8$. The sign condition is $4/y-4/z<0$, i.e.\ $y>z$, and
$\sign A=\sign(16/y-8/z)=\sign(2z-y)$: the canard point is supercritical
for $y>2z$ and subcritical for $z<y<2z$. At $(x,y,z)=(1,64,1)$ and
$(1,\frac43,1)$ the first Lyapunov coefficient confirms both, with
$\sqrt\eps\,l_1\to-0.248$ and $1.000$. This splitting arises when $\ce Y$
and $\ce Z$ are abundant. If their concentrations are of order $1/\eps$ and
the rate constant of $\ce X+\ce Y+\ce Z\to2\ce Z$ is of order $\eps^2$, this
reaction consumes $\ce X$ at a rate of order one, while the relative changes
of $\ce Y$ and $\ce Z$ that it causes are of order $\eps$; with a rate
constant of order $\eps$, $\ce Z\to\ce Y$ changes them on the same slow time
scale.
\end{example}

\begin{remark}[How small $\eps$ must be]\label{rem:BT}
The results above concern the limit $\eps\to0$. At a Hopf point
of~\eqref{eq:planar} near a canard point, $F_z+\eps\chi_z=-\eps G_u$, and the
Jacobian on the class has determinant
$-\eps\big(\eps(G_u^2+\chi_uG_z)+F_uG_z\big)$. Evaluating this at the
canard point, as in \citet[Remark~2.5]{BeldimanHalanay2026}, gives the
estimate
\[
  \eps_{\mathrm{BT}}=-\frac{F_uG_z}{G_u^2+\chi_uG_z}
\]
for the value of $\eps$ at which the singular Hopf branch ends at a double
zero eigenvalue, if the denominator is positive. The branch can also end at
a fold of the Hopf curve, and the estimate is accurate only when it is small:
in Example~\ref{ex:both} the branch ends at $0.75$ and $0.88$ times the
estimates $0.086$ and $0.016$. With one fast species,
$G_u=\sum_iU_ih_im^{(i)}_{01}$ by the proof of Proposition~\ref{prop:sup},
so for $\chi=0$, that is, when the slow parts of the reactions have no
component along $w$,
\[
  \eps_{\mathrm{BT}}=-\frac{h_{\ce X}\,m_{10}\sum_iU_ih_in^{(i)}_{01}}
  {\big(\sum_iU_ih_im^{(i)}_{01}\big)^2}.
\]
Like the criticality, how small $\eps$ must be therefore depends on the
composition of the slow pool (Example~\ref{ex:both}). The estimate is not
a bound, and we tested it numerically: at random canard
points of three-species networks the branch ends within a factor two of
$\eps_{\mathrm{BT}}$ at $123$ of the $135$ points with
$\eps_{\mathrm{BT}}<0.03$, and at fixed fluxes the end follows
$\eps_{\mathrm{BT}}$ as the ratio of the slow species changes
(Section~\ref{sm:eps}).
\end{remark}

Consider the universe of the $84$ reactions on $\ce X,\ce Y,\ce Z$ with
$\ce X$ fast, $\ce Y+\ce Z$ conserved and molecularity at most three. A
complete enumeration of its elementary modes shows that the dependence on
the slow pool is common: of the $1\,302$ modes that are canard networks on
their own, $150$ are supercritical at every
admissible point, $584$ subcritical and $8$ degenerate, while $560$ admit
both criticalities, depending on $z/y$. As an independent check,
Theorem~\ref{thm:modes} agrees with direct exact linear programming on
networks sampled from this universe (Section~\ref{sm:one}).

\section{Sign-locking}\label{sec:lock}

A canard point requires a fold point on which the slow flux vanishes. In
many networks this is impossible because the slow flux has a fixed sign on
the fold for all rate constants. With one fast species, the fold points
are the fast equilibria with singular fast Jacobian; in fluxes they form
the \emph{fold-flux cone}
$\Phi=\{v>0:\ \sum_rL_rv_r=\sum_rL_ra_rv_r=0\}$, independent of the point,
and unless $n_{20}$ vanishes identically on $\Phi$ the nondegenerate folds,
those with $n_{20}\neq0$, form an open dense subset of $\Phi$. We say
that $\beta$ is \emph{sign-locked} if the network has a nondegenerate fold
and $\beta$ has the same nonzero sign at all of them, for all rate
constants. Since $\beta$ is linear and $\Phi$ is a relatively open convex
cone, $\beta$ is positive on $\Phi$, negative on $\Phi$, identically zero on
$\Phi$, or of both signs, in which case the canard cone $K$ is nonempty.

\begin{theorem}[Certificate]\label{thm:farkas}
If $\Phi\neq\emptyset$, then $\beta\ge0$ on $\Phi$ iff there are
$\lambda_1,\lambda_2$ with $g_r\ge L_r(\lambda_1+\lambda_2a_r)$ for every
reaction, i.e.\ iff $\beta-\lambda_1F-\lambda_2x\,\partial_xF$ has
nonnegative coefficients. In that case $\beta>0$ on $\Phi$ unless all these
inequalities are equalities, and then $\beta$ vanishes identically on the
fold.
\end{theorem}

\begin{proof}
$\Phi$ is the relative interior of $C=\{v\ge0:\ Lv=(L\circ a)v=0\}$, and by
Farkas' lemma $g\cdot v\ge0$ on $C$ iff $g-\lambda_1L-\lambda_2L\circ a\ge0$
for some $\lambda$. On $\Phi$, $\beta=\sum_rc_rv_r$ with $c\ge0$ and $v>0$.
\end{proof}

\begin{figure}[ht]
\centering
\includegraphics[width=\textwidth]{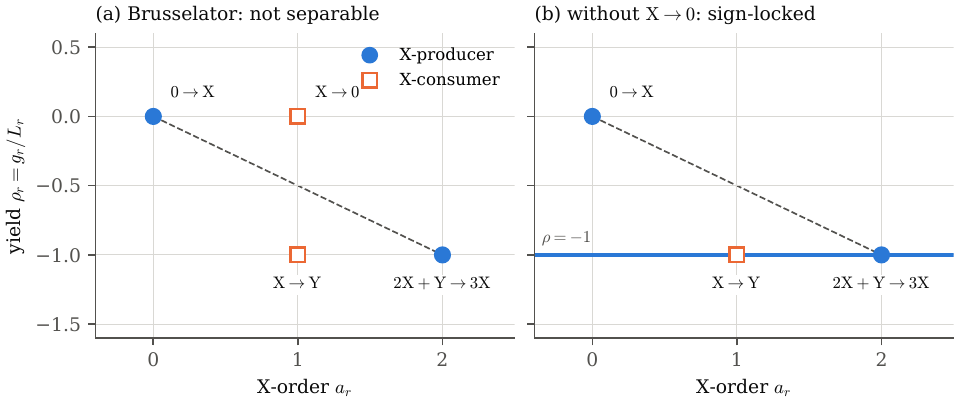}
\caption{Yield diagrams. Each reaction that changes $\ce X$ is drawn at its
$\ce X$-order $a_r$ and its slow yield $\rho_r=g_r/L_r$. (a) The
Brusselator: the consumers $\ce X\to\ce Y$ and $\ce X\to0$ lie on both sides
of the chord through the producers, so no line separates, $\beta$ changes
sign on the fold, and there is a canard point. (b) Without $\ce X\to0$ the
line $\rho=-1$ separates, and $\beta+F=k_1>0$: the slow flux is
positive at every fast equilibrium, for all rate constants.}
\label{fig:yield}
\end{figure}

For a reaction with $L_r\neq0$ let $\rho_r=g_r/L_r$ be its slow yield per
unit change of $\ce X$, and call it an $\ce X$-producer if $L_r>0$ and an
$\ce X$-consumer if $L_r<0$. By Theorem~\ref{thm:farkas}, $\beta\ge0$ on the
fold iff in the \emph{yield diagram} of the points $(a_r,\rho_r)$ a line
lies on or below the producers and on or above the consumers, and no
reaction with $L_r=0$ consumes the slow species (Figure~\ref{fig:yield}).

\begin{theorem}[Yield diagram]\label{thm:yield}
If $\Phi\neq\emptyset$, then $\beta\ge0$ on $\Phi$ iff (Y1) every reaction
with $L_r=0$ has $g_r\ge0$, (Y2) every producer has a yield at least that of
every consumer of the same order, and (Y3) in every alternating triple of
orders $a_1<a_2<a_3$, with the outer reactions of one kind and the middle
one of the other, the middle point lies on or below the chord through the
outer points if these are producers, and on or above it if they are
consumers.
\end{theorem}

\begin{proof}
$\beta\ge0$ on $C$ iff $g(e)\ge0$ for every elementary mode $e$ of
$[L;\,L\circ a]$, which have support at most three. Support one is a
reaction with $L_r=0$, giving (Y1). Support two forces a producer $p$ and a
consumer $c$ of equal order, with $g(e)=L_p|L_c|(\rho_p-\rho_c)$, giving
(Y2). Support three requires distinct orders and alternating kinds; the
mode is the cross product of the two rows, and
$g(e)=-L_1L_2L_3\det\begin{psmallmatrix}1&1&1\\a_1&a_2&a_3\\\rho_1&\rho_2&\rho_3\end{psmallmatrix}$,
which gives (Y3).
\end{proof}

\begin{corollary}[Hierarchy]\label{cor:hier}
A network with one fast species and a nondegenerate fold is exactly one of
the following:
\begin{enumerate}[label=(\roman*),nosep]
  \item sign-locked, so that no fold point is an equilibrium;
  \item with $\beta\equiv0$ on the fold, the case $\rank B=2$, in which every
        fold point is an equilibrium and canard points, if any, have
        codimension zero;
  \item with $\beta$ changing sign on the fold and a canard point, which
        happens iff the mode criterion holds;
  \item \emph{flat}, with $\beta$ changing sign on the fold but
        $n_{20}\equiv0$ on $K$, so that every fold equilibrium is degenerate;
  \item \emph{saddle-locked}, with $\beta$ changing sign on the fold,
        $n_{20}\not\equiv0$ on $K$ and $F_uG_z\ge0$ at every fold
        equilibrium, so that no fold equilibrium is a canard point; those
        with $F_uG_z>0$ are saddles.
\end{enumerate}
Cases (i), (iv) and (v) are obstructions to canard points.
The modes of a subnetwork are modes of the network, so $\beta\ge0$ on the
fold passes to subnetworks, and a sign change of $\beta$ on the fold of a
subnetwork with a nondegenerate fold passes to every supernetwork.
\end{corollary}

The yield diagram is a two-dimensional picture of a Farkas certificate, but
its content is dynamical. A separating line is a combination
$\beta-\lambda_1F-\lambda_2x\,\partial_xF$ with nonnegative coefficients, so
at every fold point, for all rate constants, the slow variable moves in the
same direction. The consequences are as follows.

\begin{proposition}[Dynamics of sign-locked networks]\label{prop:lockdyn}
Let $w=e_{\ce X}$ and suppose that $\beta$ is sign-locked. Then, for all rate
constants and on every class:
\begin{enumerate}[label=(\alph*),nosep]
  \item at every point with $F=\beta=0$ we have $F_z\neq0$: no equilibrium of
        the reduced problem lies on a fold;
  \item for every compact set $\mathcal K$ of pairs of rate constants and
        points of $\R^n_{>0}$ there is $\eps_0>0$ such that, for
        $0<\eps<\eps_0$, no equilibrium $x$ of~\eqref{eq:planar} with
        $(k,x)\in\mathcal K$ is a Hopf point; in particular there is no
        singular Hopf bifurcation and no canard explosion;
  \item along a continuous family of rate constants, a continuous branch of
        points with $F=\beta=0$ stays on branches of $C_0$ of one type,
        attracting or repelling.
\end{enumerate}
\end{proposition}

\begin{proof}
(a) A point with $F=F_z=0$ has fluxes in $\Phi$. Since $\beta$ has one
nonzero sign at the nondegenerate folds, which are dense in $\Phi$, it is
nonnegative or nonpositive on $\Phi$, and by Theorem~\ref{thm:farkas} it
vanishes nowhere on $\Phi$. (b) If not, there are $\eps_n\to0$ and Hopf points
$x_n$ of the system with rate constants $k_n$ and $\eps=\eps_n$, with
$(k_n,x_n)\in\mathcal K$. At $x_n$ we have $F+\eps_n\chi=0$, $G=0$ and, since the
trace of the Jacobian on the class vanishes, $F_z+\eps_n(\chi_z+G_u)=0$. A
subsequence converges to some $(k,x)\in\mathcal K$ with $F=G=F_z=0$, which
contradicts~(a). (c) $F_z$ is continuous along the branch and does not
vanish by~(a).
\end{proof}

So for a sign-locked network the yield diagram decides, for all rate
constants at once, that oscillations cannot be born in a singular Hopf
bifurcation near a fold, and an equilibrium on the repelling branch cannot
be moved to an attracting branch along a continuous branch of solutions of
$F=\beta=0$, such as one along which $\det\partial(F,\beta)/\partial(z,u)\neq0$.
The statement~(b) is uniform on compact sets: for each choice of rate
constants, and on compact sets of them, it holds for $\eps$ small enough.

The fold is decided by the orders and kinds of the reactions
\citep{JoshiShiu2017}, \citep[Proposition~4.1]{BeldimanHalanay2026}, and
sign-locking by configurations of at most three reactions; flatness and
saddle-locking are decided by the modes of $B$, that is, by configurations
of up to twelve reactions. As an independent check, direct linear
programming, the certificate, the modes and the yield diagram agree on
random networks of the planar universe of \citet{BeldimanHalanay2026} and of
the three-species universe of Section~\ref{sec:one}
(Section~\ref{sm:lock}). In these samples, each network consists of one
cooperative $\ce X$-producer and two to six further reactions drawn at
random, and about two thirds of them have no fold. Among those that fold,
$42$--$46\%$ are sign-locked, $31$--$33\%$ have canard points,
$20$--$25\%$ are saddle-locked, $0.7$--$2.4\%$ have $\beta\equiv0$ on the
fold and one network is flat. These frequencies depend on how the
networks are drawn and are used in no proof, but they show that, among
networks that fold, sign-locking is a common obstruction. The same
mechanism appears outside the rank-one class. Network~26 of the smallest
Hopf networks \citep{BanajiBoros2023}, $\ce X\to2\ce X$,
$2\ce X\to2\ce Y$, $\ce Y\to2\ce Z$, $\ce X+\ce Z\to0$, has with $\ce Y$
slow the fast species $\ce X$ and $\ce Z$ and a fast subnetwork of rank
two. On its critical set the fold condition reads $2v_2-v_4=0$, and
$\beta=2v_2-v_3=-\tfrac12(2v_3-v_4)+(2v_2-v_4)+\tfrac12v_4$ is a
combination of the equation of $\ce Z$ and of this condition plus $v_4/2$,
so $\beta=k_4xz/2>0$ at every fold point.

\section{The global canard explosion}\label{sec:global}

The flux criterion gives a nondegenerate singular Hopf bifurcation of known
criticality. The full explosion needs the global hypotheses of
\citet{KrupaSzmolyan2001relax} on the critical curve $C_0=\{F=0\}$ of the
class: (A1) $C_0$ is a graph $u=\varphi(z)$ with exactly two nondegenerate
critical points, a minimum and a maximum; (A2) the outer branches attract
and the middle branch repels; (A3$'$) at $\lambda=0$ one fold is a
nondegenerate canard point and the other a nondegenerate fold with $G\neq0$;
(A4$'$) the slow flow runs from the attracting branch through the canard
point and towards the other fold on the far branch. Their Theorems~3.3 and
3.5 then give a family of periodic orbits that passes from Hopf cycles
through canard cycles to relaxation oscillations within an exponentially
small $\lambda$-interval.

\begin{theorem}[Global canard explosion]\label{thm:global}
Suppose the class has a canard point with $\mathcal Q_h\neq0$, unfolded by a
parameter $\lambda$ of the rate constants. Let $D$ be the intersection of the class with the positive
orthant and $C$ the connected component of $C_0\cap D$ that contains the
canard point. Assume (G1) $F_u\neq0$ on $C$; (G2) $C$ has exactly two
folds, both nondegenerate, and a repelling middle branch; (G3) $\beta$
vanishes on $C$ only at the canard point; and (G4) each outer branch of $C$
reaches the heights of both folds, so that the singular cycles of Krupa and
Szmolyan exist, these cycles lie in $D$, and their fast fibres meet no other
component of $C_0\cap D$. Then the system on the
class undergoes the canard explosion of \citet[Theorems~3.3 and
3.5]{KrupaSzmolyan2001relax}, supercritical if $\mathcal Q_h<0$ and
subcritical if $\mathcal Q_h>0$.
\end{theorem}

\begin{proof}
By (G1) and the implicit function theorem $z$ is a local coordinate along
the connected curve $C$, which is therefore a graph $u=\varphi(z)$ over an
interval, with $\varphi'=-F_z/F_u$. Since $F_u$ has one sign on $C$,
$F_z=-F_u\varphi'$ changes sign exactly at the folds, so the outer branches
attract iff the middle one repels; with reflections $z\mapsto-z$,
$u\mapsto-u$, which preserve all conditions and the sign of $A$, the canard
point becomes the minimum, and (G1)--(G2) give (A1)--(A2). The canard point
is nondegenerate by Theorem~\ref{thm:flux}, and the other fold has
$\beta\neq0$ by (G3), so (A3$'$) holds. For (A4$'$), let $S_l$, $S_m$ and
$S_r$ be the left, middle and right branches of $C$. Near the canard point
(A2) gives $F_{zz}>0$, the minimum gives $\varphi''=-F_{zz}/F_u>0$ and hence
$F_u<0$, and the sign condition gives $G_z>0$. So $G<0$ on $S_l$ and $G>0$ on
$S_m$ near the canard point, and by (G3) these signs hold on all of $S_l$
and on $S_m\cup S_r$. With $\varphi'<0$ on $S_l$ and $S_r$ and $\varphi'>0$
on $S_m$, the reduced flow $z'=G/\varphi'$, in slow time, has the signs
required by (A4$'$); without (G3), a further zero of $G$ on $C$ would be an
equilibrium of the reduced flow. Finally, by (G4) the singular cycles
exist, and the construction of Krupa and Szmolyan uses only a compact
neighbourhood of them, which lies in $D$ and meets $C_0$ only in $C$. The conditions are imposed at $\lambda=0$; since
the folds are nondegenerate and $F_u\neq0$, they persist on this
neighbourhood for $\lambda$ near $0$, as (A1)--(A2) require.
\end{proof}

The theorem also says where the explosion happens. Let $\mu$ be the
unfolding parameter, with the canard point at $\mu=\mu_0$, and let $G_\mu$ be
the derivative of $G$ along the continued fold point.

\begin{corollary}[Location of the explosion]\label{cor:location}
Under the assumptions of Theorem~\ref{thm:global}, the singular Hopf point
$\mu_H(\eps)$ and the maximal canard $\mu_c(\eps)$ satisfy
\[
  \mu_c(\eps)-\mu_H(\eps)=S\eps+O(\eps^{3/2}),\qquad
  S=\frac{G_z\,\mathcal Q_h}{2F_{zz}^3\,G_\mu}.
\]
The explosion lies on the side of $\mu_H$ on which the Hopf cycles exist.
\end{corollary}

\begin{proof}
By \citet[Theorems~3.1 and 3.2, (3.15)--(3.16)]{KrupaSzmolyan2001relax},
$\lambda_c-\lambda_H=-\tfrac18A\eps+O(\eps^{3/2})$ in the parameter of their
canonical form. After the translation that keeps the fold at the origin, the
rescaling in the proof of Lemma~\ref{lem:A} gives
$Y'/\eps=X+(\tau/b)G_\mu(\mu-\mu_0)+\cdots$, so
$\lambda=-(\tau/b)G_\mu(\mu-\mu_0)$ up to terms of higher order and a shift of
order $\eps$ common to $\lambda_H$ and $\lambda_c$. Hence
$\mu_c-\mu_H=bA\eps/(8\tau G_\mu)+O(\eps^{3/2})$, and $A=2\tau\mathcal Q_h/F_{zz}^2$
and $b=2G_z/F_{zz}$ give $S$. In the canonical form, increasing $\lambda$
moves the equilibrium onto the repelling branch, so for $A<0$ the explosion
lies beyond the loss of stability, where the stable Hopf cycles live, and
for $A>0$ on the side of the stable equilibrium, where the unstable ones
live.
\end{proof}

So the width of the parameter interval of small cycles, and the side of the
Hopf point on which they live, are again given by forms in the fluxes.

\begin{figure}[ht]
\centering
\includegraphics[width=\textwidth]{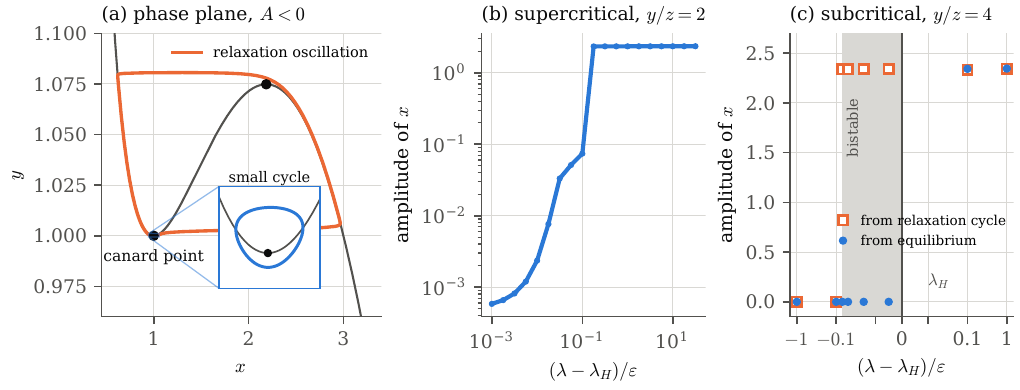}
\caption{The canard explosion of Example~\ref{ex:both} for $\eps=10^{-3}$.
(a) Supercritical case: the critical curve, its two folds, the relaxation
oscillation at $\lambda=\lambda_H+3\cdot10^{-4}$ and (inset) the small cycle at
$\lambda=\lambda_H+10^{-4}$, before the jump; the intermediate canard cycles
live in an exponentially small $\lambda$-interval and are not resolved.
(b) Amplitude of the attractor in the
supercritical case: small cycles grow and jump to the relaxation
oscillation. (c) Subcritical case: amplitude of the attractor reached from
the equilibrium and from the relaxation cycle; in the shaded window both
coexist.}
\label{fig:explosion}
\end{figure}

\begin{example}[One network, both explosions]\label{ex:both}
Let $\ce X$ be fast and $\ce Y+\ce Z=T$ conserved, with fast reactions
$0\to\ce X$, $\ce X+\ce Y\to\ce Y$, $2\ce X\to3\ce X$, $3\ce X\to2\ce X$ and slow
reactions $\ce X+\ce Z\to\ce X+\ce Y$, $\ce Y\to\ce Z$, that is, reactions
whose rate constants are of order $\eps$. Then
$F=k_1-k_2xy+k_3x^2-k_4x^3$ and $\beta=k_5x(T-y)-k_6y$, and the critical
curve is the graph $y=\varphi(x)=(k_1+k_3x^2-k_4x^3)/(k_2x)$. The fluxes
$(s-2,2s-3,s,1,m,m)$ with $s>2$ and $m>0$ lie in $K$ at $x=1$; for
$s\neq3$, so that $F_{zz}=2(s-3)\neq0$, the folds are at $x=1$ and
$x=\big(s-2+\sqrt{(s-2)(s+6)}\big)/4$, and by
Proposition~\ref{prop:crit} $\sign A=\sign\big((2s-6)^2y-18(s-2)z\big)$. For
$s=5$ the folds are at $1$ and $(3+\sqrt{33})/4$, and $A<0$ iff $y/z<27/8$.
The fast reactions satisfy conditions (S1)--(S4) below but not (S6), which
fails for $\ce Z$, so by Theorem~\ref{thm:folds} conditions (G1), (G2)
and (G4) hold whenever the other fold lies in the class, below $y=T$. For
$s=5$ and $y=1$ this means $y/z<(413+77\sqrt{33})/64\approx13.4$. Condition
(S5) fails, and (G3) has to be checked at given rate constants. At
$(x,y,z)=(1,1,\frac12)$ and $(1,1,\frac14)$, i.e.\ $y/z=2$ and $4$, the
numerator of $\beta$ on the critical curve is $(x-1)(2x^3-7x^2+9x+3)$ and
$(x-1)(4x^3-15x^2+15x+3)$, and exact real-root isolation shows that the cubic
factors have no positive root, so $\beta$ vanishes on the critical curve
only at the canard point. Figure~\ref{fig:explosion}
shows the explosion for $\eps=10^{-3}$, with $k_6\mapsto k_6(1+\lambda)$. At
$y/z=2$ the attracting cycle grows to amplitude $0.074$ and becomes a
relaxation oscillation of amplitude $2.34$ within
$\lambda-\lambda_H\in(1.0,1.8)\cdot10^{-4}$. At $y/z=4$ the stable
equilibrium coexists with the relaxation oscillation for
$\lambda_H-7\cdot10^{-5}\lesssim\lambda<\lambda_H$, and beyond $\lambda_H$
the four orbits we integrated (two values of $\lambda$, two initial
conditions each) go to the relaxation oscillation. Here $m=1$, which gives
$A=-1.04$ and $A=0.47$; in general $A$ is proportional to $\sqrt m$. For
$s=5$ and $y=1$ one finds $\mathcal Q_h=m(16/z-54)$, $G_z=m$, $F_{zz}=4$ and,
for the unfolding $k_6\mapsto k_6(1+\lambda)$, $G_\lambda=-m$, so
Corollary~\ref{cor:location} gives $\lambda_c-\lambda_H\approx m(54-16/z)\,\eps/128$, that is, for $m=1$,
$1.7\cdot10^{-4}$ beyond the Hopf point at $y/z=2$ and $7.8\cdot10^{-5}$
before it at $y/z=4$, for $\eps=10^{-3}$, in agreement with the simulations.
So the composition of the slow pool decides not only whether oscillations
set in softly or with hysteresis, but also on which side of the Hopf point,
and over which parameter interval, the small cycles exist.
At $x=1$ the estimate of Remark~\ref{rem:BT} is
$\eps_{\mathrm{BT}}=(2s-3)yz^2/\big(m(y+z)^2\big)$, which is $7/(1+y/z)^2$ at
$y=1$ for $s=5$ and $m=1$; so the slow pool also sets how small $\eps$ must
be. At $y/z=2$ and $4$ the singular Hopf branch, continued in $\eps$,
persists up to $\eps=0.23$ and $0.14$, where it turns back in $\eps$. At
$y/z=8$ and $20$ it ends at a double zero eigenvalue at $\eps=0.064$ and
$0.014$, against the estimates $0.086$ and $0.016$.
Condition (G3) is a genuine restriction: for $s=4$ the critical curve is
shallow, with folds at $1$ and $1.618$, and at $y/z=20$ the slow flux has
two further zeros on the middle branch, at $x\approx1.407$ and $1.587$.
\end{example}

For networks with one fast species whose fast equation is affine in the
slow variable, three of the four conditions follow from sign conditions on
the reactions as soon as the second fold lies in the class, and (G3) follows
from injectivity on the class. Let $w=e_{\ce X}$ and $U_{\ce X}=0$ as in
Section~\ref{sec:one}, let $a_r$ be the $\ce X$-order of reaction $r$, and
suppose that every reaction with $L_r\neq0$ has total order at most one in
the slow species. On a class the slow species have concentrations
$x_i=c_i+U_iu$, and species that are neither fast nor slow have constant
concentrations, which we absorb into the rate constants. Each reaction with
$L_r\neq0$ then has a factor $\ell_r(u)$, equal to $x_i$ if it has order one
in the slow species $i$ and to $1$ otherwise, and
\[
  F=\mathcal A(x)+u\,\mathcal B(x),
\]
with polynomials $\mathcal A$ and $\mathcal B$ in $x=x_{\ce X}$. The
critical curve is the graph of $\varphi=-\mathcal A/\mathcal B$ wherever
$\mathcal B\neq0$, and its folds are the positive roots of
\[
  P=\mathcal A'\mathcal B-\mathcal A\mathcal B'
  =\sum_{r<s}k_rk_sL_rL_s(a_r-a_s)D_{rs}\,x^{a_r+a_s-1},
\]
summed over pairs of reactions with $L_r,L_s\neq0$, where
$D_{rs}=\ell_r\ell_s'-\ell_s\ell_r'$ is constant. Thus $D_{rs}=0$ if $r$
and $s$ involve the same slow species or none, $D_{rs}=U_j$ if only $s$
involves one, $j$, $D_{rs}=-U_i$ if only $r$ involves one, $i$, and
$D_{rs}=x_iU_j-x_jU_i$, which has the sign of $U_j$
when $U_iU_j<0$, if they involve different slow species $i$ and $j$.

At given rate constants, the hypotheses of Theorem~\ref{thm:global} are a
finite computation in this setting.

\begin{proposition}[Deciding (G1)--(G4)]\label{prop:decide}
Let $w=e_{\ce X}$, and suppose that every reaction with $L_r\neq0$ has total
order at most one in the slow species. For given rational rate constants and
a given class, (G1)--(G4) at a canard point are decided by exact real root
isolation of univariate polynomials with rational coefficients and by
comparison of real algebraic numbers.
\end{proposition}

\begin{proof}
On the class, $D$ is the product of $x>0$ and an interval $I$ of $u$. Where
$\mathcal B(x)=0$, either $\mathcal A(x)\neq0$ and $C_0$ has no point, or
$C_0$ contains the whole line, on which $F_u=\mathcal B=0$; elsewhere $C_0$
is the graph of $\varphi=-\mathcal A/\mathcal B$. So the components of
$C_0\cap D$ on which $F_u\neq0$ are the pieces of this graph over the
intervals of $x$ bounded by the positive roots of $\mathcal B$ and of
$\mathcal A+u_e\mathcal B$, $u_e$ an endpoint of $I$, and (G1) holds iff the
component $C$ of the canard point is such a piece and no line of $C_0$ in
$D$ meets it. Its folds are the roots of $P$ in its interval, nondegenerate
iff simple, and the branch types follow from the sign of
$F_x=-\mathcal B\varphi'$, which gives (G2). The zeros of $\beta$ on $C$ are
the roots of the numerator of $\beta(x,\varphi(x))$ in the interval, which
gives (G3). For (G4), the range of $\varphi$ on each outer branch of $C$ is
determined by its value at the fold and its limit at the endpoint, both real
algebraic numbers, and the singular cycles exist iff the interval between the
fold values lies in both ranges. They then lie in $D$, since $C$ does, and
their fast fibres are horizontal segments between two points of $C$, so they
stay over the interval of $x$ of $C$ and meet no other component of
$C_0\cap D$, all of which lie over disjoint sets of $x$.
\end{proof}

Without the affine structure the same conditions are first-order statements
about polynomials and remain decidable, by cylindrical algebraic
decomposition, but no longer by univariate computations. Example~\ref{ex:both}
is an instance of Proposition~\ref{prop:decide}, and so is the two-variable
Oregonator (Proposition~\ref{prop:oreg2} below).

A coefficient of $P$, or of another polynomial below, is a polynomial in the
rate constants and the concentrations, and its sign is \emph{determined} if
all its terms have the same sign for all rate constants and classes. We use
the following conditions.
\begin{enumerate}[label=(S\arabic*),nosep]
  \item Every reaction with $L_r\neq0$ has total order at most one in the
        slow species, and the nonzero numbers $L_rU_iy_{ri}$ have one sign
        $\sigma$.
  \item The coefficients of $P$ have determined signs with exactly two sign
        changes.
  \item Let $p_\pm$ be the lowest and highest $\ce X$-orders of the
        reactions with $L_r\neq0$ of order zero in the slow species, and
        $q_\pm$ those of order one. Then $p_-<q_-$ and $p_+>q_+$, the
        reactions of the first kind of order $p_-$ have $L_r$ of one sign,
        and so do those of order $p_+$, and the infinite limits
        $\varphi(0^+)$, $\varphi(\infty)$ have opposite signs.
  \item $\sigma\sign\varphi(0^+)<0$.
  \item $U_i(L_rg_s-L_sg_r)(a_ry_{si}-a_sy_{ri})\ge0$ for every slow
        species $i$ and every pair of reactions, with at least one positive.
  \item For every slow species $j$, the reactions with $L_r\neq0$ of order
        zero in $j$ have, at each $\ce X$-order, $L_r$ of one sign, and these
        signs change at most once as the order increases.
\end{enumerate}
For $n=2$, with $\ce Y$ slow and $U=e_{\ce Y}$, $\mathcal A$ and
$\mathcal B$ collect the reactions of $\ce Y$-order zero and one, (S6) says
that the reactions in $\mathcal A$ have, order by order, one sign of $L_r$,
changing at most once, and (S5) is the
classical injectivity condition $(L_rg_s-L_sg_r)(a_rb_s-a_sb_r)\ge0$, with
$b_r$ the $\ce Y$-order \citep{CraciunFeinberg2005,ShinarFeinberg2012}. For
more species, (S5) is that condition on the class, summed over the slow
species as in Proposition~\ref{prop:sup}; we include the short proof
because it gives (G3).

\begin{theorem}[One equilibrium]\label{thm:inj}
Let $w=e_{\ce X}$ and assume (S5). Then on every class $F=\beta=0$ has at
most one solution in $D$, for all rate constants; in particular (G3) holds
at every canard point.
\end{theorem}

\begin{proof}
Let $f=(F,\beta)$ and let $x\neq x'$ be points of $D$. For every $r$,
$x^{y_r}-x'^{y_r}=c_r\ip{y_r}{\ln x-\ln x'}$ with $c_r>0$, and
$\ln x-\ln x'=(x-x')/\xi$ componentwise, with $\xi$ the logarithmic mean of
$x$ and $x'$. Since $x-x'=\Delta z\,e_{\ce X}+\Delta u\,U$, this gives
$f(x)-f(x')=M(\Delta z,\Delta u)^{\mathsf T}$, where
\[
  M=\sum_r\tilde v_r\begin{pmatrix}L_r\\ g_r\end{pmatrix}
  \Big(\frac{a_r}{\xi_{\ce X}},\ \sum_i\frac{U_iy_{ri}}{\xi_i}\Big),
  \qquad \tilde v_r=k_rc_r>0,
\]
as in the standard injectivity argument
\citep{CraciunFeinberg2005,ShinarFeinberg2012}. By the Cauchy--Binet
formula,
\[
  \det M=\frac1{\xi_{\ce X}}\sum_i\frac1{\xi_i}\sum_{r<s}\tilde v_r\tilde v_s\,
  U_i(L_rg_s-L_sg_r)(a_ry_{si}-a_sy_{ri})>0,
\]
so $f(x)\neq f(x')$.
\end{proof}

\begin{theorem}[Two folds, shape and room]\label{thm:folds}
Let $w=e_{\ce X}$ and assume (S1)--(S4). Then (G1), (G2) and (G4) hold at
every canard point whose other fold lies in $D$. Under (S6) it always does.
\end{theorem}

\begin{proof}
By (S1), $\mathcal B$ has all coefficients of sign $\sigma$, so $\varphi$ is
defined on $x>0$ and $F_u=\mathcal B\neq0$, which is (G1). Since
$P=-\mathcal B^2\varphi'$, the folds are the positive roots of $P$; by
Descartes' rule there are at most two, counted with multiplicity, and their
number is even, so the simple root at the canard point forces exactly two
simple roots. By (S3), $\varphi$ runs between the infinite limits and the
fold values on the outer branches, so it has a minimum and a maximum, in an
order fixed by (S3), and since $F_x=-\mathcal B\varphi'$ on the curve, (S4)
makes the branch next to $x=0$, and with it the other outer branch,
attracting. In the coordinates $(x,u)$, $D$ is the product of $x>0$ and an
open interval $I$ of $u$. If both fold values lie in $I$, each of the three
monotone branches takes each endpoint value of $I$ at most once, so
$C_0\cap D$ is a single arc that contains both folds and on which both outer
branches reach every height between the fold values; the fast fibres
$u=\mathrm{const}$ are horizontal segments in $D$. This gives (G2) and
(G4). Finally, at a finite endpoint $u_e$ of $I$ some slow species vanish
and the others are positive, so $F(\cdot,u_e)$ is a polynomial in $x$ whose
terms come from reactions of order zero in a vanishing species $j$; if
several species vanish, from a subset of these, whose signs order by order
form a subsequence. By (S6) its coefficients change sign at most once, so $\varphi=u_e$ has at most
one positive solution, counted with multiplicity. If the other fold had its
value at $u_e$ or beyond, while the canard fold has its value in $I$, then
$\varphi-u_e$ would vanish on all three branches, or doubly at the other
fold; so the other fold lies in $D$.
\end{proof}

\begin{corollary}\label{cor:structural}
Let $w=e_{\ce X}$ and $\rank B=3$. Under (S1)--(S6), every canard point of
the network with $\mathcal Q_h\neq0$, for any rate constants and on any
class, unfolds into the full canard explosion.
\end{corollary}

The smallest networks with these properties have four reactions. For two
species one is $0\to\ce X$, $\ce X\to\ce Y$, $2\ce X+\ce Y\to3\ce X$,
$3\ce X\to\ce Y$, for which $\mathcal A=k_1-k_2x-3k_4x^3$,
$\mathcal B=k_3x^2$ and $P=k_3x(-2k_1+k_2x-3k_4x^3)$; for three species,
with $\ce Y+\ce Z$ conserved, one is $0\to3\ce X$,
$2\ce X+\ce Y\to3\ce X+\ce Z$, $\ce X+\ce Z\to\ce Y$, $4\ce X\to0$.
A complete enumeration of the $138\,722$ minimal canard networks of
codimension one of \citet{BeldimanHalanay2026} finds $63$ that satisfy all
the conditions: $30$
supercritical, $5$ subcritical, $27$ with both criticalities and one with
$\mathcal Q_h\equiv0$ on its canard cone, to which the corollary does not
apply. Conditions that hold for all rate constants and on all classes are
necessarily rare, because for most networks the shape of the critical
curve depends on the rate constants; the bottleneck is the S-shape, (S2), and
not the other conditions (Section~\ref{sm:struct}). For given rate
constants Theorem~\ref{thm:global} with Proposition~\ref{prop:decide}
applies without these restrictions, as in Example~\ref{ex:both} and
Proposition~\ref{prop:oreg2}. For three species, (S5) is restrictive. In the universe of
Section~\ref{sec:one} the only cooperative producer, $2\ce X\to3\ce X$,
violates it together with every reaction that changes $\ce Y$, so (G3) has
to be checked at given rate constants, as in Example~\ref{ex:both}; with
catalysed producers such as $2\ce X+\ce Y\to3\ce X+\ce Z$ it can hold.
The hypothesis on the other fold in Theorem~\ref{thm:folds} is needed: among
random canard points of three-species networks with (S1)--(S4), it fails at
$45$, and at each of them (G2) and (G4) fail, while at the other $891$ (G1),
(G2) and (G4) hold, as the theorem asserts. These and further checks are
reported in Section~\ref{sm:struct}.

The structural conditions can give both explosions in one network, for all
rate constants.

\begin{example}[Both explosions, for all rate constants]\label{ex:structboth}
Let $\ce X$ be fast and $\ce Y$ slow in
$0\to\ce X$, $\ce X\to\ce Y$, $2\ce X+\ce Y\to3\ce X$, $3\ce X\to0$,
$3\ce X\to3\ce Y$, the Brusselator with its linear removal of $\ce X$
replaced by two cubic ones. The network satisfies (S1)--(S6) and has
$\rank B=3$. Its canard cone is
$K=\{(3v_4,\,9v_4+3v_5,\,9v_4+6v_5,\,v_4,\,v_5):v_4,v_5>0\}$, on which
$n_{20}=-6v_5$, $n_{01}m_{10}=-27v_4(3v_4+2v_5)<0$ and
$Q=-162\,v_4(3v_4-v_5)(3v_4+2v_5)$, in the moment forms of
Section~\ref{sec:flux} for which $\mathcal Q_h=Q/(x^4y)$. So every flux vector of
$K$ gives a canard point, and since $v_5/v_4=k_5/k_4$ at every point,
Corollary~\ref{cor:structural} shows that every canard point of the network,
for any rate constants with $k_5\neq3k_4$, unfolds into the full canard
explosion, which is supercritical if $k_5<3k_4$ and subcritical if
$k_5>3k_4$. The first
Lyapunov coefficient confirms both cases (Section~\ref{sm:struct}).
\end{example}

Here the criticality is set by a ratio of rate constants; in
Example~\ref{ex:both} it is set by the composition of the slow pool, but that
network violates (S5) and (S6). Among networks on three species that satisfy
all of (S1)--(S6) we found no composition switch: at $11\,640$ canard flux
vectors of $366$ such networks, $Q^{(\ce Y)}$ and $Q^{(\ce Z)}$ have opposite
signs or one of them vanishes, so that by Proposition~\ref{prop:crit} the
criticality does not depend on the composition (Section~\ref{sm:struct}).
Whether the structural conditions exclude the switch is open.

We close the section with a classical model.

\begin{proposition}[The two-variable Oregonator]\label{prop:oreg2}
Let $\ce X$ be fast and $\ce Y$ slow in the network
\begin{gather*}
  \ce X\to2\ce X,\ 2\ce X\to3\ce X,\ 3\ce X\to2\ce X,\ \ce X+\ce Y\to\ce Y,\
  \ce Y\to\ce X+\ce Y,\\
  \ce X\to\ce X+\ce Y,\ 2\ce X\to2\ce X+\ce Y,\ \ce X+\ce Y\to\ce X,\ \ce Y\to0,
\end{gather*}
with rate constants $q,1-q,1,f,fq$ for the fast reactions in the first line
and $q,1,1,q$ for the slow ones in the second, so that
$F=qx+(1-q)x^2-x^3-fxy+fqy$ and $\beta=(x-y)(x+q)$. Let $0<q<q^*$, where
$q^*=0.07973\ldots$ is the real root of $10q^3+6q^2+12q-1$. Then the critical
curve has exactly two folds in the positive quadrant, at $x_-<x_+$, both
nondegenerate; for each of them exactly one $f>0$ makes it a canard point;
$f$ unfolds it; and every such canard point with $\mathcal Q_h\neq0$
undergoes the full canard explosion. The canard point at $x_-$ is
supercritical, and the one at $x_+$ is subcritical for $q<q_c$ and
supercritical for $q_c<q<q^*$, where $q_c=0.05551\ldots$ is the only root in
$(0,q^*)$ of $16q^6+152q^5+68q^4+190q^3+25q^2-20q+1$.
\end{proposition}

\begin{proof}
The critical curve is $y=\varphi(x)=x(1-x)(x+q)/\big(f(x-q)\big)$, since
$F=2q^2(1-q)\neq0$ on the line $x=q$. As $\varphi<0$ for $0<x<q$ and for
$x>1$, $C_0\cap D$ is the graph of $\varphi$ over $(q,1)$, on which
$F_y=-f(x-q)\neq0$: this is (G1), and $C_0\cap D=C$. Next,
$\varphi'=P/\big(f(x-q)^2\big)$ with $P=-2x^3+(1+2q)x^2-2q(1-q)x-q^2$, and
$P'(q)=0$, $P(q)=-2q^2(1-q)<0$, $P(1)=q^2-1<0$, while $P$ has a local maximum
at $x=(1-q)/3$ with value $-(10q^3+6q^2+12q-1)/27>0$ for $q<q^*$. So $P$ has
exactly two roots in $(q,1)$, both simple: $\varphi$ decreases from $+\infty$
to a minimum at $x_-$, increases to a maximum at $x_+$, and decreases to $0$
at $x=1$. Since $F_x=-F_y\varphi'=f(x-q)\varphi'$ on the curve, the outer
branches attract and the middle one repels, which is (G2), and the outer
branches cover $(\varphi(x_-),\infty)$ and $(0,\varphi(x_+))$, which with
$C_0\cap D=C$ gives (G4). On the curve, $\beta=0$ iff $y=x$, that is,
$(1-x)(x+q)=f(x-q)$, a quadratic whose roots have product $-q(1+f)<0$; so
$\beta$ vanishes on $C$ at exactly one point, which gives (G3), and the fold
$x_\pm$ is a canard point iff $f=(1-x_\pm)(x_\pm+q)/(x_\pm-q)>0$. There
$F_yG_x=-(x+q)^2(1-x)<0$. The folds do not depend on $f$, and
$\varphi\propto1/f$, so the derivative of $\beta$ along the fold point with
respect to $f$ is $(x_\pm+q)\varphi(x_\pm)/f\neq0$. At the canard point with
fold abscissa $x$,
\[
  \mathcal Q_h=\frac{2(x+q)\big(12x^3-6qx^2-3x^2-12q^2x+7qx-x-2q^3+6q^2-q\big)}{x-q}
\]
and its numerator
has a common root with $P$ iff the resultant
$-64q^4(q^2-1)(16q^6+152q^5+68q^4+190q^3+25q^2-20q+1)$ vanishes, that is,
in $(0,q^*)$, iff $q=q_c$. So the sign of $\mathcal Q_h$ at each canard
point is constant on $(0,q_c)$ and on $(q_c,q^*)$, and its values at
$q=1/40$ and $q=7/100$ give the statement.
\end{proof}

Multiplying the vector field by $1/(x+q)>0$ changes time but not the orbits,
so the same holds for the reduction of \citet{TysonFife1980} of the
Oregonator, in which bromide is eliminated. As $q\to0$ the canard points
tend to $f=1+\sqrt2$ and $f=1/2$, the classical ends of the oscillatory
range. The first Lyapunov coefficient confirms all signs at
$q=1/40,1/20,6/100,7/100$ (Section~\ref{sm:struct}).

\section{A fast direction that mixes species}\label{sec:mixing}

Let $\mathcal J=\{i:w_i\neq0\}$ and $n^{(i)}(v)=\sum_rL_ry_{ri}v_r$. Then
$c^d_r=\sum_{i\in\mathcal J}y_{ri}w_ih_i$, so the fold form is
\[
  F_z=\sum_{i\in\mathcal J}w_ih_in^{(i)}(v),
\]
which is linear in $v$ for a fixed point and linear in $(h_i)_{i\in\mathcal J}$ for
fixed fluxes, while the rest of the canard cone,
$K_0=\{v>0:\ Lv=0,\ gv=0\}$, does not depend on the point. For two mixing
species, $\mathcal J=\{1,2\}$, put $p=w_1n^{(1)}$, $q=w_2n^{(2)}$ and
$t=h_1/h_2=x_2/x_1$; the fold condition is $tp+q=0$. A \emph{fold
equilibrium} is a point with $F=F_z=\beta=0$.

\begin{theorem}[Fold ratio]\label{thm:ratio}
Suppose $K_0\neq\emptyset$. Fold equilibria with $x_2/x_1=t$ exist, for
some rate constants, iff the line $tp+q=0$ meets the cone
\[
  \Sigma=\{(p(v),q(v)):v\in K_0\}=\relint\operatorname{cone}\{(p(e),q(e))\},
\]
where $e$ runs over the elementary modes of $[L;\,g]$. The admissible ratios form the empty set, a point, all
of $(0,\infty)$, or at most two open intervals whose finite endpoints are
slopes $-q(e)/p(e)$ of mode images.
\end{theorem}

\begin{proof}
The image of the relatively open cone $K_0$ is the relatively open cone
$\Sigma$. If $0\in\Sigma$, every ratio is admissible. Otherwise its directions form an open arc of length at most $\pi$
or a single direction, and the line meets it in the open quadrants II and
IV in at most two open intervals of $t$, bounded by mode images.
\end{proof}

For networks A and B of Section~\ref{sec:flux}, $K_0$ gives
$t=(v_1+v_3)/(v_3-v_1)$ and $t=(v_2+v_5+v_6)/(v_6-v_2)$, so in both
canard points require $z>x$, as observed on the grid. At a fixed point all
conditions of Theorem~\ref{thm:flux}(a) are linear in $v$, which gives the
mode criterion for every rank-one fast subnetwork.

\begin{corollary}[Modes at a point]\label{cor:pointmodes}
Assume (H1)--(H2) and $K_h\neq\emptyset$, $h=1/x$. A network has a canard
point at $x$ iff the elementary modes of $B_h$ include a mode with
$F_{zz}\neq0$ and a possibly equal pair $e,e'$
for which the cone spanned by $(F_u(e),G_z(e))$ and $(F_u(e'),G_z(e'))$
meets $\Omega=\{\varphi\psi<0\}$. Every network with a canard point
contains a subnetwork with a canard point carried by at most three modes of
some $B_h$, hence with at most twelve reactions.
\end{corollary}

\begin{proof}
By Theorem~\ref{thm:flux}(a), $x$ is a canard point for some rate constants
iff some $v\in K_h$ satisfies $F_{zz}(v)\neq0$ and $F_u(v)G_z(v)<0$, with
all forms evaluated at $h$. These are the conditions of the planar mode
criterion with $n_{20}$, $n_{01}$ and $m_{10}$ replaced by the linear forms
$F_{zz}$, $F_u$ and $G_z$ at $h$, and the proof of Theorem~\ref{thm:modes},
with $p_i$ replaced by $(F_u,G_z)$, applies to the cone $K_h$; its
elementary modes have support at most four, since $B_h$ has three rows.
\end{proof}

With one fast species $B_h$ has the kernel and the modes of $B$, and
Corollary~\ref{cor:pointmodes}, taken over all points, reduces to
Theorem~\ref{thm:modes}. For two mixing species $B_h$ depends on $h$ only
through the fold ratio $t$, up to scaling a row, with entries affine in $t$,
so the supports of its modes change only at finitely many ratios; the sign
condition still depends on the concentrations of the slow species.

Sign-locking carries over piecewise. Over all points, the fluxes of fold
points are those with $Lv=0$ and $tp+q=0$ for some $t>0$, that is, the
union of the three relatively open convex cones
\begin{gather*}
  P_+=\{v>0:Lv=0,\ p>0>q\},\qquad P_-=\{v>0:Lv=0,\ p<0<q\},\\
  P_0=\{v>0:Lv=0,\ p=q=0\},
\end{gather*}
where we write $p$ and $q$ also for their coefficient vectors. We call the
network \emph{locked} if the slow flux vanishes at no fold point, for any
point and any rate constants; unlike in Section~\ref{sec:lock}, its sign may
differ between the pieces.

\begin{proposition}[Piecewise certificates]\label{prop:mixlock}
The network is locked iff on each nonempty piece the slow flux has a
certificate of one sign: for $P_+$, a sign $s\in\{\pm1\}$, a number $\lambda$
and $\mu,\nu\ge0$ such that $sg-\lambda L-\mu p+\nu q\ge0$ componentwise and
$(sg-\lambda L-\mu p+\nu q,\mu,\nu)\neq0$; for $P_-$ the same with $\mu p$
and $\nu q$ of opposite signs; and for $P_0$ the same with $\mu,\nu\in\R$
and the last condition replaced by $sg-\lambda L-\mu p-\nu q\neq0$.
\end{proposition}

\begin{proof}
On a nonempty relatively open convex cone the linear form $\beta=g\cdot v$
has no zero iff it is positive or negative there, and $s\beta>0$ on $P_+$
implies $s\beta\ge0$ on the closure
$\{v\ge0:Lv=0,\ p\ge0\ge q\}$. By Farkas' lemma this holds iff
$sg=\lambda L+\mu p-\nu q+c$ with $\mu,\nu\ge0$ and $c\ge0$. Then
$s\beta=c\cdot v+\mu p(v)-\nu q(v)$ on $P_+$, which is positive unless $c$,
$\mu$ and $\nu$ all vanish. The pieces $P_-$ and $P_0$ are treated in the
same way.
\end{proof}

For example, $\ce Z\to\ce X$, $\ce X\to\ce Y$, $2\ce X+\ce Y\to3\ce X$ with
$\ce X+\ce Y+\ce Z$ conserved is locked, with $\beta+F=k_1z>0$, and adding
$\ce X\to\ce Z$ removes the locking.

As an independent check, Theorem~\ref{thm:ratio},
Corollary~\ref{cor:pointmodes} and Proposition~\ref{prop:mixlock} agree with
exact linear programming on
random networks of a universe of $126$ reactions on $\ce X,\ce Y,\ce Z$
with a mixing fast direction (Section~\ref{sm:mixing}).

\section{Fast subnetworks of higher rank}\label{sec:rank}

We now drop (H1) and keep one slow direction. Suppose that the fast parts of
the reactions span a subspace of dimension $s\ge1$, with basis the columns of
$W\in\R^{n\times s}$, so that~\eqref{eq:sys} becomes
\begin{equation}\label{eq:sysW}
  \dot x=W\sum_rL_rv_r(x)+\eps\sum_rN^1_rv_r(x),\qquad L_r\in\R^s,
\end{equation}
and assume
\begin{enumerate}[label=(H2$'$),nosep]
  \item $\mathcal S=\operatorname{im}W+\operatorname{span}\{N^1_r\}$ has
        dimension $s+1$.
\end{enumerate}
With $U$, $\varpi$ and $g_r$ as in Section~\ref{sec:setting}, now with
$\varpi W=0$, the coordinates $x=x^*+Wy+uU$ turn~\eqref{eq:sysW} on the class
through $x^*$ into
\begin{equation}\label{eq:planarW}
  y'=f(y,u)+\eps\chi(y,u),\qquad u'=\eps\beta(y,u),\qquad
  f=\sum_rL_rv_r,\quad \beta=\sum_rg_rv_r,
\end{equation}
with $y\in\R^s$. The critical set on the class is the curve $f=0$, and its
fold points are those at which the fast Jacobian $D_yf$ is singular. For
$s=1$ this is~\eqref{eq:planar}. The Oregonator below has $s=2$.

\paragraph{Flux forms.}
As in Lemma~\ref{lem:deriv}, for directions $a,b\in\R^n$ the derivatives of
$v_r$ at $x=1/h$ are $D v_r[a]=v_rc_r[a]$, $D^2v_r[a,b]=v_rc_r[a,b]$ and
$D^3v_r[a,a,a]=v_rc_r[a,a,a]$, with
\begin{gather*}
  c_r[a]=\ip{y_r}{a\circ h},\qquad
  c_r[a,b]=c_r[a]c_r[b]-\ip{y_r}{a\circ b\circ h^2},\\
  c_r[a,a,a]=c_r[a]^3-3c_r[a]\ip{y_r}{a^2\circ h^2}+2\ip{y_r}{a^3\circ h^3}.
\end{gather*}
For weights $c\in\R^R$ and fluxes $v$ write $L[c]=\sum_rL_rc_rv_r\in\R^s$ and
$g[c]=\sum_rg_rc_rv_r$. The fast Jacobian on the class is
\begin{equation}\label{eq:M}
  M(v,h)=\big(L[c[W_1]],\dots,L[c[W_s]]\big)
  =\sum_rv_r\,L_r\,(y_r\circ h)^{\mathsf T}W,
\end{equation}
where $W_j$ is the $j$th column of $W$. Suppose that $0$ is an algebraically
simple eigenvalue of $M$, let $\xi$ and $\ell$ be right and left null vectors
with $\ell\xi=1$, and let $M^\sharp=(M+\xi\ell)^{-1}-\xi\ell$ be the group inverse,
which inverts $M$ on its range $\{\ell\eta=0\}$ and annihilates $\xi$. Put
$d=W\xi$, the fold direction in species space,
\[
  \Psi_u=-M^\sharp L[c[U]],\qquad \Psi_{zz}=-M^\sharp L[c[d,d]],\qquad
  \hat U=U+W\Psi_u,\qquad \omega=W\Psi_{zz},
\]
and define
\begin{equation}\label{eq:formsW}
\begin{gathered}
  F_{zz}=\ell L[c[d,d]],\quad F_u=\ell L[c[U]],\quad F_{zu}=\ell L[c[d,\hat U]],\quad
  F_{zzz}=\ell L[c[d,d,d]]+3\,\ell L[c[d,\omega]],\\
  G_z=g[c[d]],\qquad G_u=g[c[\hat U]],\qquad G_{zz}=g[c[d,d]]+g[c[\omega]],
\end{gathered}
\end{equation}
and $\mathcal Q$ by~\eqref{eq:Qh} with these forms. For $s=1$, $W=w$, we
have $M=F_z$ and, at a fold, $\xi=\ell=1$ and $M^\sharp=0$, so $\Psi_u=\Psi_{zz}=0$
and~\eqref{eq:formsW} are the forms of Section~\ref{sec:flux}. In general the
forms differ from those of Section~\ref{sec:flux}, taken along $d$, only
through $\Psi_u$ and $\Psi_{zz}$, which are the first- and second-order
slaving of the transverse fast directions to the fold direction. By this we
mean the following: near the fold the transverse fast directions do not
evolve independently, but on the centre manifold constructed in the proof of
Theorem~\ref{thm:rank} they are given as a graph
$y=y^*+z\xi+\Psi(z,u)$, with $\ell\Psi=0$, over the fold coordinate $z$ and the
slow coordinate $u$, and $\Psi_u$ and $\Psi_{zz}$ are the derivatives of this
graph that enter $A$.

A \emph{canard point} of~\eqref{eq:planarW} is a point with $f=0$ and
$\beta=0$ at which (C1) $0$ is an algebraically simple eigenvalue of
$M$ and no other eigenvalue lies on the imaginary axis, and (C2)
$F_{zz}\neq0$ and $F_uG_z<0$, with the forms evaluated at the fluxes of the
point. It is \emph{unfolded} by a parameter $\lambda$ of the rate constants if
the derivative of $\beta$ along the continued fold point does not vanish.

\begin{theorem}[Flux criterion for one slow direction]\label{thm:rank}
Assume that the fast parts of the reactions span $\operatorname{im}W$ and
that (H2$'$) holds.
\begin{enumerate}[label=(\alph*),nosep]
  \item For given rate constants, $x$ is a canard point iff its fluxes $v$
        satisfy $Lv=0$, $gv=0$, (C1) and (C2) at $h=1/x$. Some rate constants
        give a canard point at $x$ iff these conditions have a solution
        $v>0$.
  \item Let $x^*$ be a canard point at $\lambda=0$, unfolded by $\lambda$,
        with $\mathcal Q\neq0$. Near $x^*$, the system on the class has a
        locally invariant manifold of dimension two for each small
        $(\lambda,\eps)$, which contains every solution that stays near $x^*$
        and attracts nearby solutions if the nonzero eigenvalues of $M$ have
        negative real parts. The flow on it is a planar slow--fast system with
        a canard point in the sense of \citet[(3.2)]{KrupaSzmolyan2001relax}
        and coefficient
        \[
          A=\frac{2\mathcal Q}{F_{zz}^2\sqrt{-F_uG_z}} .
        \]
        For every small $\eps>0$ the system on the class has a unique
        equilibrium near $x^*$, which undergoes a Hopf bifurcation at some
        $\lambda_H(\sqrt\eps)\to0$, supercritical if $\mathcal Q<0$ and
        subcritical if $\mathcal Q>0$.
  \item The null vectors can be chosen rational in $(v,h)$ and homogeneous of
        degree zero in $v$; then the forms~\eqref{eq:formsW} are rational in
        $(v,h)$ and homogeneous of degree one in $v$. The fold condition
        $\det M(v,h)=0$ is homogeneous of degree $s$ in $v$. The sign of
        $\mathcal Q$ and the value of $A$ do not depend on the normalisation
        of $\xi$ and $\ell$.
\end{enumerate}
\end{theorem}

\begin{proof}
(a) Every $v>0$ is realised at every $x$ by $k_r=v_r/x^{y_r}$. By the chain
rule, $M=D_yf$, $L[c[U]]=D_uf$, $L[c[a,b]]=D^2f[a,b]$ and so on, with
derivatives in species space along the indicated directions, so (C1) and
(C2) are conditions on $(v,h)$.

(b) Append $\lambda'=0$ and $\eps'=0$. At $(x^*,0,0)$ the linearisation has
the eigenvalues of $M$ and three further zeros, and since $\beta(x^*)=0$ its
generalised kernel is the four-dimensional space of
$(y-y^*,u-u^*,\lambda,\eps)$ with $y-y^*=z\xi+\eta$, $\ell\eta=0$ and
$\eta=-M^\sharp(D_uf\,(u-u^*)+f_\lambda\lambda+\chi\eps)$; by (C1) the other
eigenvalues are off the imaginary axis. By the centre manifold theorem
\citep[Ch.~5]{Kuznetsov2004} there is, for any $k\ge3$, a locally invariant
$C^k$ manifold $W^c=\{y=y^*+z\xi+\Psi(z,u-u^*,\lambda,\eps)\}$ with $\ell\Psi=0$,
tangent to this space, which contains every solution that stays near the
point for all $t\in\R$, is
attracting under the stated condition, and is foliated by the invariant
slices of constant $(\lambda,\eps)$. In particular it contains the layer
equilibria near $x^*$ and the fold point. On $W^c$ the coordinate
$z=\ell(y-y^*)$ satisfies
\[
  z'=\ell f(y^*+z\xi+\Psi,u)+\eps\,\ell\chi=:\hat F,\qquad
  u'=\eps\,\beta(y^*+z\xi+\Psi,u)=:\eps\hat G,
\]
a planar system of the form (3.1) of \citet{KrupaSzmolyan2001relax}. At
$\lambda=\eps=0$, invariance of $W^c$ projected by $P=I-\xi\ell$ reads
$Pf(y^*+z\xi+\Psi_0,u)=\partial_z\Psi_0\,\hat F$, where $\Psi_0=\Psi(\cdot,\cdot,0,0)$.
Since $\hat F=O(z^2)+O(u-u^*)$, comparing the terms of first order in $u$ and of
second order in $z$ gives $M\Psi_u=-PD_uf$ and $M\Psi_{zz}=-PD^2f[d,d]$ for the
derivatives of $\Psi_0$, which with $\ell\Psi=0$ are the expressions above. Differentiating
$\hat F$ and $\hat G$ by the chain rule, every term that contains a further
derivative of $\Psi_0$, such as $\Psi_{zu}$ or $\Psi_{zzz}$, is multiplied by
$\ell M=0$, and the remaining terms are the forms~\eqref{eq:formsW}; also
$\hat F_z=\ell M\xi=0$ and $\hat F_u=\ell D_uf=F_u$. So the reduced system has
a fold at the canard point with these coefficients. The Jacobian of
$(\hat F,\hat F_z)$ in $(z,u)$ there has determinant $-F_uF_{zz}\neq0$, so the
fold point continues smoothly in $\lambda$; it lies in $W^c$, so the
derivative of $\hat G$ along it is the derivative of $\beta$ along the
continued fold point of the system on the class. After the
$\lambda$-dependent translation that keeps the fold at the origin, which
changes no coefficient at $\lambda=\eps=0$, (C2) and the unfolding are the
nondegeneracy conditions (3.2) of Krupa and Szmolyan, and $A$
follows as in Theorem~\ref{thm:flux}(b); the term $\eps\ell\chi$ enters only
$a_1$. By their Theorem~3.1 the reduced system has a Hopf bifurcation at
$\lambda_H(\sqrt\eps)$ of the criticality of $A$. Equilibria of the full
system near $x^*$ lie in $W^c$. At such an equilibrium the tangent plane of the
slice carries the two eigenvalues of the reduced system, and the other $s-1$
eigenvalues are close to the nonzero eigenvalues of $M$. So at
$\lambda_H$ the slice is a centre manifold of the Hopf point, and the Hopf
bifurcation of the system on the class has the same criticality
\citep[Ch.~5]{Kuznetsov2004}.

(c) By (C1) the adjugate of $M$ has rank one, $\operatorname{adj}M=\kappa\,\xi_0\ell_0$
with $\ell_0\xi_0=1$ and $\kappa=\operatorname{tr}\operatorname{adj}M\neq0$, the
product of the nonzero eigenvalues of $M$. If $(\operatorname{adj}M)_{ij}\neq0$,
then $\xi=(\operatorname{adj}M)_{\cdot j}/\kappa$ and
$\ell=(\operatorname{adj}M)_{i\cdot}/(\operatorname{adj}M)_{ij}$ are null
vectors with $\ell\xi=1$, rational in $(v,h)$ and, since $M$ is linear in
$v$, homogeneous of degree zero. Then $M^\sharp$ is homogeneous of degree
$-1$, $\Psi_u$ and $\Psi_{zz}$ of degree $0$, and each form of degree one.
Replacing $(\xi,\ell)$ by $(c\xi,\ell/c)$ multiplies $F_{zz}$ and $G_z$ by
$c$, $F_u$ by $1/c$, $F_{zzz}$ and $G_{zz}$ by $c^2$, leaves $F_{zu}$ and $G_u$
unchanged, and so multiplies $\mathcal Q$ by $c^2$ and leaves $A$ unchanged.
\end{proof}

By the Cauchy--Binet formula,
\[
  \det M(v,h)=\sum_{|S|=s}\det L_S\,\det\big((y_r\circ h)^{\mathsf T}W\big)_{r\in S}
  \prod_{r\in S}v_r,
\]
summed over sets $S$ of $s$ reactions, with $L_S$ the columns of $L$ in $S$.
So the canard fluxes at a point are cut out of the polyhedral cone
$\{v>0:Lv=0,\ gv=0\}$ by a form of degree $s$, whose terms come from sets of
$s$ reactions. This is why the polyhedral theory of
Sections~\ref{sec:one}--\ref{sec:mixing}, with its elementary modes and
Farkas certificates, is a theory of rank one: for $s=1$ the fold condition is
linear in the fluxes, while for $s\ge2$ it has in general degree $s$; it can
factor, as for Network~26 of Section~\ref{sec:lock}, where on $Lv=0$ it
reduces to $2v_2-v_4=0$. What carries over to every rank is
the reduction of the local problem to forms in $(v,h)$ that are homogeneous in
$v$, and the realisation argument of part~(a). The corrections also have a
simple structure. The direct cubic term $\ell L[c[d,d,d]]$ vanishes unless
some reactant has molecularity at least three, since $c_r[a,a,a]$ is the third
derivative of a monomial of degree $|y_r|$. The transverse term
$3\ell L[c[d,\omega]]=-3\,\ell D^2f\big[d,WM^\sharp D^2f[d,d]\big]$ couples
pairs of reactions through $M^\sharp$, that is, through the part of the fast
subnetwork transverse to the fold. In the supplementary materials
(Section~\ref{sm:rank}) we show that Theorem~\ref{thm:rank}(b) agrees with
the intrinsic criterion of \citet{DeMaesschalckDoanWynen2021} applied on the
centre manifold.

Theorem~\ref{thm:rank} is local. The global argument of
Theorem~\ref{thm:global} uses that the system on the class is planar. For
$s\ge2$ it would need, in addition, a two-dimensional attracting invariant
manifold along the singular cycles, which one may expect from Fenichel
theory when the transverse fast directions contract uniformly along the
critical curve and the fast fibres; we do not pursue this here.

\begin{example}[The Oregonator]\label{ex:oregonator}
The Oregonator \citep{FieldNoyes1974} in the scaling of
\citet[eqs.~(14)--(16)]{BronsBarEli1991} is
\[
  \eps'\dot x=x+y-xy-2qx^2,\qquad \dot y=fz-y-xy,\qquad \dot z=\eps(x-z),
\]
where $x$, $y$ and $z$ are the scaled concentrations of HBrO$_2$, bromide
and the oxidised catalyst, $q=1/40$, and $\eps'$ is the ratio of the time
scales of HBrO$_2$ and bromide. It is the mass-action system of the network
\[
  \ce Y\to\ce X,\quad \ce X+\ce Y\to0,\quad \ce X\to2\ce X+\ce Z,\quad
  2\ce X\to0,\quad \ce Z\to2\ce Y,\quad \ce Z\to0
\]
with rate constants $1,1,1,q,f/2,1-f/2$, for $0<f<2$, in which the changes
of $\ce X$ and $\ce Y$ are fast and those of $\ce Z$ slow; as in
Example~\ref{ex:onemode}, the slow parts come from an abundant species. So
$s=2$, $W=(e_{\ce X},e_{\ce Y})$ with the $\ce X$-row of $L$ divided by $\eps'$,
$U=e_{\ce Z}$, and there is no conservation law. Write $M$ and $\ell$ for
the Jacobian~\eqref{eq:M} and its left null vector at $\eps'=1$. Dividing the
$\ce X$-row of $L$ by $\eps'$ replaces $L$ by $DL$ and $M$ by $DM$, with
$D=\operatorname{diag}(1/\eps',1)$. This changes neither the kernel of $L$
nor the zero set of $\det M$, so the canard points do not depend on $\eps'$.
Theorem~\ref{thm:rank}(a) gives exactly two, at $f=0.7135870120$ and
$f=1.8532510835$: the point $\mathrm M$ of \citet[Table~II]{BronsBarEli1991},
at $(x,y,z)=(8.6777,0.6399,8.6777)$, and a point $\mathrm N$ at
$(2.7425,1.3581,2.7425)$, whose $(x,y)$ is their fold point $\mathrm N$ and
which is a canard point at this value of $f$. The matrix $DM$ has the right
null vector $\xi$ of $M$ and a left null vector proportional to $\ell D^{-1}$,
and normalising it gives $\ell_{\eps'}DL=c(\eps')\,\ell L$ with
$c(\eps')=\operatorname{tr}M/(M_{11}+\eps'M_{22})$. At both points
$M_{11}<0$ and $M_{22}=-1-x<0$ ($M_{11}=-0.508$ and $-0.632$), so
$c(\eps')>0$ and $\operatorname{tr}DM\neq0$, and (C1) and (C2) hold for every
$\eps'$. If the corrections through $M^\sharp$ were dropped, that is
$\Psi_u=\Psi_{zz}=0$, every $F$-form would be multiplied by $c(\eps')$ and
every $G$-form left unchanged, so $A$ would be divided by $\sqrt{c(\eps')}$ and
keep its sign. So $\eps'$ can change the sign of $A$ only through $M^\sharp$.
Table~\ref{tab:oregonator} gives $A$ and the value it would have without the
corrections.

\begin{table}[ht]
\centering
\caption{The coefficient $A$ of Theorem~\ref{thm:rank} at the canard points
$\mathrm M$ and $\mathrm N$ of the Oregonator, and its value without the
corrections through $M^\sharp$ ($\Psi_u=\Psi_{zz}=0$). The sign of the first
Lyapunov coefficient of the three-variable system agrees with that of $A$ in
all $24$ cases: six values of $\eps'$, two points, and $\eps=10^{-3}$ and
$10^{-4}$.}\label{tab:oregonator}
\small
\begin{tabular}{@{}lrrrrrr@{}}
\toprule
$\eps'$ & $1$ & $\tfrac12$ & $\tfrac15$ & $\tfrac1{10}$ & $\tfrac1{100}$ & $\tfrac1{1000}$\\
\midrule
$A$ at $\mathrm M$ & $0.457$ & $0.110$ & $-0.246$ & $-0.413$ & $-0.317$ & $-0.183$\\
without corrections & $2.727$ & $1.976$ & $1.336$ & $1.038$ & $0.664$ & $0.615$\\
$A$ at $\mathrm N$ & $-6.757$ & $-4.646$ & $-2.563$ & $-1.365$ & $0.514$ & $0.796$\\
without corrections & $2.328$ & $1.761$ & $1.308$ & $1.117$ & $0.911$ & $0.888$\\
\bottomrule
\end{tabular}
\end{table}

Both canard points change criticality as $\eps'$ varies, at
$\eps'=0.388375$ for $\mathrm M$ and $\eps'=0.028855$ for $\mathrm N$; the first Lyapunov
coefficient of the three-variable system changes sign at $0.388411$ and
$0.028854$ for $\eps=10^{-5}$. Without the corrections, both points would be
subcritical for every $\eps'$. The limit $\eps'\to0$ is the planar reduction of
\citet{BronsBarEli1991}, in which HBrO$_2$ is slaved, that is, eliminated
through its quasi-steady state $\eps'\dot x=0$, and there $\mathrm M$ is
supercritical and $\mathrm N$ subcritical, as they found. In the three-variable
model, $\eps'=1$, it is the other way round. The opposite limit
$\eps'\to\infty$, in which bromide is slaved, that is, eliminated through
$\dot y=0$, is the two-variable Oregonator
of Proposition~\ref{prop:oreg2} with $q=1/20$, after the scaling
$x\mapsto x/20$, $z\mapsto z/20$; its canard points are at the same values of
$f$, and there, as for $\eps'=1$ and for $\eps'=10$ and $100$ (where $A=2.42$ and
$7.98$ at $\mathrm M$, and $-21.5$ and $-68.0$ at $\mathrm N$), $\mathrm M$ is
subcritical and $\mathrm N$ supercritical. For $\eps'=1$ and
$\eps=10^{-3}$, direct simulation shows the subcritical explosion at $\mathrm M$: below the Hopf point
the relaxation oscillation coexists with the stable equilibrium over an
interval of length $9.4\cdot10^{-5}$ in $f$. So the criticality of a canard
explosion depends on the full fast subnetwork and not only on a reduction of
it. Theorem~\ref{thm:rank} locates the effect in the slaving of the
transverse fast direction to the fold, which enters only through $M^\sharp$:
the corrections make $\mathrm M$ supercritical for $\eps'<0.388375$ and
$\mathrm N$ supercritical for $\eps'>0.028855$.
\end{example}

\section{Discussion}\label{sec:disc}

For networks with a rank-one fast subnetwork and one slow direction, the
existence, the criticality and the global unfolding of canard points are
governed by linear forms in the fluxes. With one fast species the forms of
the planar theory reappear, summed over the slow species, and with a
mixing fast direction the point enters only through one bilinear fold
condition. In both cases minimal canard networks have at most three
elementary modes and twelve reactions. What is new is that the criticality
depends on the composition of the slow pool: according to the ratio of its
slow species, the same network can have a supercritical or a subcritical
singular Hopf bifurcation (Example~\ref{ex:onemode}) and can explode in
either way (Example~\ref{ex:both}). Without the assumption of rank one, the
local theory survives in the same form: the coefficients are forms in the
fluxes, homogeneous of degree one, with corrections through the group
inverse of the fast Jacobian (Theorem~\ref{thm:rank}). These corrections are
not small in practice. In the Oregonator they are responsible for the change
of criticality of both canard points between the planar reduction and the
three-variable model (Example~\ref{ex:oregonator}), and among random networks
of rank two they change the sign of $A$ at $8$ of the $30$ canard points we
sampled (Section~\ref{sm:rank}).

Several questions remain. With three or more mixing species the fold
condition has several free ratios, and the window of
Theorem~\ref{thm:ratio} becomes a semi-algebraic region. For fast subnetworks of
rank $s\ge2$ two problems are open. The global explosion requires a
two-dimensional attracting invariant manifold along the singular cycles, which
Theorem~\ref{thm:rank} provides only near the canard point. And the
existence of canard points is no longer decided by elementary modes, since
the fold condition is in general a form of degree $s$ in the fluxes; its terms come
from sets of $s$ reactions by the Cauchy--Binet formula, and a combinatorial
theory would have to describe how the signs of these terms meet the
polyhedral cone $\{v>0:Lv=0,\ gv=0\}$. The structural
conditions of Section~\ref{sec:global} are sufficient but not necessary,
and they assume one fast species and a fast equation affine in the slow
variable. With several slow species, injectivity on the class is
restrictive, and (G3) usually has to be checked at given rate constants, as
in Example~\ref{ex:both}, which Proposition~\ref{prop:decide} makes a finite
computation; whether the structural conditions exclude a dependence of the
criticality on the composition of the slow pool, as our samples suggest, is
open. Finally, a complete
classification of a three-species universe would require an enumeration of
the modes of $B_h$ along the fold ratio, whose combinatorial type is
piecewise constant with finitely many breakpoints.

\paragraph{Data and code availability.}
The Julia and Python scripts that produce all computations and figures of
this paper, the logs of all runs and the data of Figure~\ref{fig:explosion}
are openly available on Zenodo \citep{CanardRankoneData2026}, at
\url{https://doi.org/10.5281/zenodo.23962730}. The README of the record
lists, for each statement checked in the supplementary materials, the script
and its log.

\paragraph{Use of AI tools.}
The authors used the AI assistant Claude (Anthropic), through Claude Code,
throughout this work: to write and run the Julia and Python code for the
computations and figures, to search the literature and check the references
against Crossref and zbMATH, to draft and check the text, including the
proofs, and to review the manuscript critically. The authors directed the
research, checked all statements, proofs and computations, and take full
responsibility for the content of the paper.

\appendix
\section{The coefficient \texorpdfstring{$A$}{A} of a planar canard point}\label{app:A}

For completeness we derive the coefficient $A$ of
\citet{KrupaSzmolyan2001relax} for a planar system
\begin{equation}\label{eq:KSplanar}
  z'=F(z,u,\lambda)+\eps\chi(z,u,\lambda,\eps),\qquad
  u'=\eps G(z,u,\lambda,\eps),
\end{equation}
with a canard point at the origin for $\lambda=\eps=0$: $F=F_z=G=0$,
$F_{zz}\neq0$ and $F_uG_z<0$, all derivatives taken at
$(z,u,\lambda,\eps)=0$, unfolded by $\lambda$: the derivative $G_\lambda$ of $G$
along the continued fold point does not vanish. Their canonical form~(3.3) is
\[
  x'=-yh_1+x^2h_2+\eps h_3,\qquad y'=\eps(xh_4-\lambda h_5+yh_6),
\]
with $h_1,h_2,h_4,h_5$ equal to $1$ and $h_3$ equal to $0$ at the origin, and
their coefficient is $A=-a_2+3a_3-2a_4-2a_5$, with $a_2=\partial_xh_1$,
$a_3=\partial_xh_2$, $a_4=\partial_xh_4$ and $a_5=h_6$ at the origin.

\begin{lemma}\label{lem:A}
At a canard point of~\eqref{eq:KSplanar},
\[
  A=\frac{2\mathcal Q}{F_{zz}^2\sqrt{-F_uG_z}},\qquad
  \mathcal Q=(G_zF_{zu}-G_uF_{zz})F_{zz}+F_u(G_{zz}F_{zz}-G_zF_{zzz}).
\]
\end{lemma}

\begin{proof}
Put $z=\alpha X$, $u=bY$ and $t=\tau T$ with
\[
  b=\frac{2G_z}{F_{zz}},\qquad \alpha^2=-\frac{4F_uG_z}{F_{zz}^2}>0,\qquad
  \tau=\frac2{\alpha F_{zz}},
\]
choosing the sign of $\alpha$ so that $\tau>0$; time is not reversed, so
the criticality is preserved. Then $\tau bF_u/\alpha=-1$,
$\tau\alpha F_{zz}/2=1$ and $\tau\alpha G_z/b=1$, so at $\lambda=\eps=0$
\[
  X'=-Y+X^2+\tau bF_{zu}XY+\tfrac16\tau\alpha^2F_{zzz}X^3+\cdots,\qquad
  \frac{Y'}{\eps}=X+\tau G_uY+\frac{\tau\alpha^2G_{zz}}{2b}X^2+\cdots,
\]
where the dots stand for terms that do not contribute to $a_2,\dots,a_5$:
in $X'$, terms in $Y^2$ and terms of order three other than $X^3$; in
$Y'/\eps$, terms in $XY$, $Y^2$ and of order three. Hence $a_2=-\tau bF_{zu}$,
$a_3=\tfrac16\tau\alpha^2F_{zzz}$, $a_4=\tau\alpha^2G_{zz}/(2b)$ and
$a_5=\tau G_u$. The remaining normalisations of the canonical form, the
$\lambda$-dependent translation that keeps the fold at the origin, a
reparametrisation of $\lambda$, which brings the $\lambda$-term to the form
$-\lambda h_5$ with $h_5(0)=1$ because $G_\lambda\neq0$, and shifts of $Y$ and
$\lambda$ of order $\eps$, change the vector field only where $\lambda\neq0$ or
$\eps\neq0$, so they do not change $a_2,\dots,a_5$; the term $\eps\chi$ enters
only $h_3$.
Substituting,
\[
  A=\tau\Big(bF_{zu}+\tfrac12\alpha^2F_{zzz}-\frac{\alpha^2G_{zz}}b-2G_u\Big)
   =\frac{2\tau\,\mathcal Q}{F_{zz}^2},
\]
and $\tau=2/|\alpha F_{zz}|=1/\sqrt{-F_uG_z}$.
\end{proof}

One checks from the formula that $A$ is unchanged under the reflections
$z\mapsto-z$ and $u\mapsto-u$, which reduce any canard point to the
orientation of the canonical form, and $A$ agrees with the intrinsic criterion of \citet{DeMaesschalckDoanWynen2021}
(Section~\ref{sm:rank}).

\clearpage
\setcounter{section}{0}
\renewcommand{\thesection}{SM\arabic{section}}
\setcounter{table}{0}
\renewcommand{\thetable}{SM\arabic{table}}
\begin{center}
\large\bfseries Supplementary materials
\end{center}
\addcontentsline{toc}{section}{Supplementary materials}
These materials report the computations that test the statements of the
paper. Unless stated otherwise, linear programs, elementary modes and
polynomial identities are computed in exact rational arithmetic (Oscar and
polymake), and the global conditions and continuations in floating point.
Every computation is a script of the Zenodo record of the paper
\citep{CanardRankoneData2026}, which also
holds its complete log (Table~\ref{tab:scripts}); random samples are seeded,
and the numbers below are those of the logs.

\section{The flux criterion}\label{sm:flux}

For $n=2$, $w=e_{\ce X}$ and $U=e_{\ce Y}$, the identities $F_z=n_{10}/x$,
$F_{zz}=n_{20}/x^2$, $F_u=n_{01}/y$, $G_z=m_{10}/x$, \dots,
$\mathcal Q_h=Q/(x^4y)$ between the forms of Section~\ref{sec:flux} and the
moment forms of the planar paper hold exactly on $200$ random flux vectors
over the $90$ reactions of the planar classification. Networks~A and~B of
Section~\ref{sec:flux} have canard points at $72$ and $20$ of the $216$
points of the grid $\{1,\dots,6\}^3$, with both criticalities. At one
canard point of each criticality in each network, $A$ from
Theorem~\ref{thm:flux}(b) agrees to $75$ digits with the value obtained by
exact symbolic differentiation of the vector field on the class, which
tests Lemma~\ref{lem:deriv}, and its sign agrees with that of the first
Lyapunov coefficient $l_1$ of the Hopf bifurcation on the class, which
tests part~(b); $\sqrt\eps\,l_1$ converges as $\eps$ decreases from
$10^{-3}$ to $10^{-5}$.

\section{One fast species}\label{sm:one}

The universe of Section~\ref{sec:one} consists of the $84$ reactions between
complexes of molecularity at most three on $\ce X,\ce Y,\ce Z$ that conserve
$\ce Y+\ce Z$; the change of $\ce X$ is fast and the change of $\ce Y$ slow,
with $U=e_{\ce Y}-e_{\ce Z}$. Its matrix $B=[L;\,L\circ a;\,g]$ has $6\,221$
elementary modes, $261$ of support two, $1\,048$ of support three and
$4\,912$ of support four. The decomposition formula of
Proposition~\ref{prop:sup} holds exactly on $200$ random pairs $(v,h)$. On
$2\,517$ sampled networks with $K\neq\emptyset$, Theorem~\ref{thm:modes}
agrees with a direct exact linear-programming decision that uses no modes:
$1\,578$ are canard networks and $939$ are not. Half of the trials draw
random sets of reactions and half draw unions of one to three elementary
modes, one of which carries a fold, so that canard networks are well
represented. Of the $1\,302$ elementary modes that are canard networks on
their own, $150$ are supercritical at every admissible point, $584$
subcritical and $8$ degenerate, while $560$ admit both criticalities,
depending on $z/y$.

\section{Sign-locking}\label{sm:lock}

On random networks of the planar universe of the planar paper and of the
three-species universe of Section~\ref{sm:one}, direct linear programming,
the certificate of Theorem~\ref{thm:farkas}, the elementary modes and the
yield diagram of Theorem~\ref{thm:yield} agree in all $1\,107$ and
$1\,150$ sampled networks with $\Phi\neq\emptyset$. Each sampled network
consists of one cooperative $\ce X$-producer ($L_r>0$, $a_r\ge2$) and two to
six further reactions drawn at random. In samples of $3\,000$ such networks
in each universe about two thirds have no fold. Among those that fold,
$42$--$46\%$ are sign-locked, $31$--$33\%$ have canard points, $20$--$25\%$
are saddle-locked, $0.7$--$2.4\%$ have $\beta\equiv0$ on the fold, and one
network is flat.

\section{The structural conditions}\label{sm:struct}

\paragraph{Two species.} On all $138\,722$ minimal canard networks of
codimension one of the planar paper, conditions (S1)--(S6) of
Section~\ref{sec:global}, evaluated by their general rules with one slow
species, agree with the planar formulation (orders $a_r$, $b_r$ and changes
of $\ce X$ and $\ce Y$). Of these networks, $34\,652$ have a fast equation
affine in $y$, $15\,389$ satisfy (S1), $150$ also (S2), $144$ also (S3) and
(S4), $141$ also (S6), and $63$ all six conditions; $543$ satisfy (S5)
alone. At random canard points, $(G1)$--$(G4)$ hold in all $126$ cases tested
for the $63$ networks. For $150$ networks with (S1) but not (S2), the
critical curve has exactly two folds in only $2$ of $300$ cases, and for
$78$ networks with (S1)--(S4) and (S6) but not (S5) the slow flux has
further zeros on the critical curve in $24$ of $156$ cases. A weaker form of
(S2), which admits coefficients of undetermined sign as long as the extreme
coefficients are determined and of one sign and no choice of the
undetermined signs gives more than three sign changes, is passed by $264$
networks instead of $150$, but none of the additional ones satisfies (S3),
(S4) and (S6); the shape of the critical curve, not the bound on the number
of folds, limits the reach of Corollary~\ref{cor:structural}.

\paragraph{Three species, molecularity three.} In the universe of
Section~\ref{sm:one} the only cooperative $\ce X$-producer is
$2\ce X\to3\ce X$, and it violates (S5) together with each of the $54$
reactions that change $\ce Y$: such a reaction consumes $\ce Y$ or $\ce Z$,
which then occurs in its reactant. So no network of this universe with a
fold satisfies (S5). Among $294\,236$ random networks, $2\ce X\to3\ce X$
plus three to seven random reactions, $400$ distinct ones satisfy
(S1)--(S4) and $\rank B=3$, and $39$ of them have canard points at the
sampled points. At their $936$ canard points with $\mathcal Q_h\neq0$, the
other fold lies in the class at $891$, and (G1), (G2) and (G4) hold at all
of them, as Theorem~\ref{thm:folds} asserts; at the other $45$, (G2) and
(G4) fail. None of these networks satisfies (S6) for both slow species.
Without an injectivity condition, (G3) holds at $652$ of the $936$ points.

\paragraph{Three species, molecularity four.} Here the cooperative
producers $2\ce X+\ce Y\to3\ce X+\ce Z$ and $2\ce X+\ce Z\to3\ce X+\ce Y$ are
compatible with (S5) together with $41$ reactions that change $\ce Y$, and
$2\ce X+\ce Y\to3\ce X+\ce Y$ and $2\ce X+\ce Z\to3\ce X+\ce Z$ with $11$;
$2\ce X\to3\ce X$, $2\ce X\to4\ce X$ and $3\ce X\to4\ce X$ with none. We drew
random sets of reactions compatible with (S5) and (S1) around a cooperative
producer. For $300$ distinct networks with (S1) and (S5), at five random
choices of rate constants and class each, the slow flux has at most one
zero on the critical curve in the class in all $1\,500$ cases, as
Theorem~\ref{thm:inj} asserts. Of $150$ networks with (S1)--(S6) and
$\rank B=3$, $47$ have canard points at the sampled points, and at their
$732$ canard points (G1), (G2), (G3) and (G4) hold at all of them.
The smallest of these networks have four reactions, for instance
$0\to3\ce X$, $2\ce X+\ce Y\to3\ce X+\ce Z$, $\ce X+\ce Z\to\ce Y$,
$4\ce X\to0$.

\paragraph{Both explosions under the structural conditions.} The network of
Example~\ref{ex:structboth} is one of the $27$ networks with both
criticalities among the $63$ that satisfy (S1)--(S6); its canard cone and
cubic were computed exactly. At $x=y=1$ with $v_4=1$ and $v_5=1$ or $5$, that
is $k_5/k_4=1$ or $5$, unfolded by $k_2$, the first Lyapunov coefficient of the
Hopf point is negative and positive, respectively, for $\eps=10^{-4}$ and
$10^{-5}$, with $\sqrt\eps\,l_1$ converging. For three species we drew
$1\,407\,969$ random sets of reactions of molecularity at most four compatible
with (S1) and (S5) around a cooperative producer, as above, and kept the
first $1\,000$ distinct networks with all of (S1)--(S6) and $\rank B=3$; $366$
of them have canard fluxes at the sampled points. At all $11\,640$ sampled
canard flux vectors, $Q^{(\ce Y)}$ and $Q^{(\ce Z)}$ have opposite signs or one
of them vanishes, so that, $U$ being $e_{\ce Y}-e_{\ce Z}$, the sign of
$h_{\ce Y}Q^{(\ce Y)}-h_{\ce Z}Q^{(\ce Z)}$ does not depend on the
composition, and no network admits both criticalities at one flux vector.

\paragraph{The two-variable Oregonator.}
The exact computations in the proof of Proposition~\ref{prop:oreg2} (the fold
polynomial, its discriminant $8q^2(q-1)(10q^3+6q^2+12q-1)$, the cubic at the
canard points and its resultant with the fold polynomial) are reproduced by a
script in exact arithmetic, which also evaluates the cubic at both canard
points for $q=1/40$, $1/20$, $6/100$ and $7/100$. The first Lyapunov
coefficient at the Hopf points near both canard points, unfolded by $f$, has
the sign of the cubic in all $16$ cases (four values of $q$, two points,
$\eps=10^{-4}$ and $10^{-5}$), with $\sqrt\eps\,l_1$ converging; in
particular the canard point at $x_+$ is subcritical for $q=1/40$ and $1/20$
and supercritical for $q=6/100$ and $7/100$, on either side of
$q_c=0.0555$.

\section{Finite \texorpdfstring{$\eps$}{epsilon}}\label{sm:eps}

For the network of Example~\ref{ex:both} with $s=5$, $m=1$ and $y=1$, the
singular Hopf branch continued in $\eps$ ends at a fold of the Hopf curve
for $y/z\le4$, at $\eps=0.46$, $0.35$, $0.23$, $0.16$ and $0.14$ for
$y/z=\frac12,1,2,\frac{27}8,4$, and at a double zero eigenvalue for
$y/z=8$ and $20$, at $\eps=0.064$ and $0.014$, against the estimates
$0.086$ and $0.016$ of Remark~\ref{rem:BT}.

For the three-species universe of Section~\ref{sm:one} we drew random
networks, $2\ce X\to3\ce X$ plus three to seven random reactions, and
random canard points with $\mathcal Q_h\neq0$, and continued the Hopf point
on the class in $\eps$ from $10^{-6}$ to $10$, with exact derivatives,
unfolded by the rate constant with the largest effect on $\beta$ along the
continued fold. At $301$ canard points of
$160$ networks the branch ends at a double zero eigenvalue at $250$, at a
fold of the Hopf curve at $41$, and persists up to $\eps=10$ at $10$. For
the $135$ points with $\eps_{\mathrm{BT}}<0.03$ the end lies within a factor
two of $\eps_{\mathrm{BT}}$ at $123$, with median ratio $0.97$; of the other
$12$, ten end earlier, seven of them at a fold of the Hopf curve and three
at a double zero eigenvalue, and two later, at $2.1$ and $2.4$ times
$\eps_{\mathrm{BT}}$. For
$0.03\le\eps_{\mathrm{BT}}<0.3$ the end is within a factor two at $87$ of
$114$ points (median ratio $0.85$), and for larger $\eps_{\mathrm{BT}}$ at
$21$ of the $42$ branches that end. At $79$ of the $301$ points
$\eps_{\mathrm{BT}}<10^{-2}$.

To isolate the slow pool we kept the fluxes of these canard points and moved
the point to $x=1$, $y+z=2$, $y/z\in\{\frac18,\frac14,\dots,8\}$; since
the canard cone does not depend on the point, $v$ stays in it, and the sign
condition decides which ratios are admissible. For $120$ flux vectors with at
least three admissible ratios, $\eps_{\mathrm{BT}}$ varies over the ratios by
a median factor of $5.3$, by more than $10$ for $34$ and by more than $100$
for $14$. For $89$ of them at least two branches with
$\eps_{\mathrm{BT}}<0.1$ end, and for $75$ of these the end lies within a
factor two of $\eps_{\mathrm{BT}}$ at every such ratio: the end of the
branch follows the composition of the slow pool. For instance, for the
network $2\ce X\to3\ce X$, $2\ce X+\ce Z\to2\ce X+\ce Y$,
$\ce Y+\ce Z\to\ce X+\ce Y+\ce Z$, $\ce X+\ce Y+\ce Z\to\ce Y+\ce Z$,
$\ce X+\ce Y\to\ce Z$, $2\ce X+\ce Z\to\ce Z$ the estimates
$0.15$, $0.22$, $0.26$, $0.23$, $0.15$, $0.075$, $0.033$ at
$y/z=\frac18,\frac14,\dots,8$ compare with ends at $0.14$, $0.18$, $0.19$,
$0.16$, $0.10$, $0.051$, $0.023$.

\section{A fast direction that mixes species}\label{sm:mixing}

The universe of Section~\ref{sec:mixing} consists of the $126$ reactions on
$\ce X,\ce Y,\ce Z$ with $\ce X+\ce Y+\ce Z$ conserved in which every change
of $\ce X$ is fast, carried by $w=e_{\ce X}-e_{\ce Z}$, and the rest is slow;
networks A and B use another splitting and do not belong to it. For $172$
random networks the set of Theorem~\ref{thm:ratio}, the union of the
intervals between consecutive mode slopes that are admissible at their
midpoints and of the admissible slopes, agrees with exact linear
programming at $174$ ratios, which confirms that admissibility changes only
at mode slopes; the set is empty for $90$ networks, connected for $81$ and
has two components for one. The pointwise mode criterion of
Corollary~\ref{cor:pointmodes} agrees with a direct linear program at $705$
random pairs of a network and a point with $K_h\neq\emptyset$, $550$ of them
with a canard point. The certificates of Proposition~\ref{prop:mixlock}
agree with the definition, checked at all critical ratios, in all $748$
pairs of a network and a sign: $148$ for random networks, none of them
locked, and $600$ for supersets of the locked example of
Section~\ref{sec:mixing}, $126$ of them locked.

\section{Fast subnetworks of higher rank}\label{sm:rank}

\paragraph{The flux forms.}
The forms~\eqref{eq:formsW} were implemented from the fluxes and $h$ alone
and compared with an implementation of the same coefficients from exact
derivatives of the vector field~\eqref{eq:planarW}, which computes $\xi$, $\ell$
and the inverse of $M$ on its range by a bordered system. At the canard
points $\mathrm M$ and $\mathrm N$ of the Oregonator, for the six values of $\eps'$ of
Table~\ref{tab:oregonator}, the two values of $A$ agree to $5\cdot10^{-15}$,
the accuracy of the double-precision null vectors of the second
implementation. Here the network is written with the integer step
$\ce Z\to2\ce Y$ of rate constant $f/2$, whereas the second implementation
uses the vector field directly.

\paragraph{The Oregonator.}
The first Lyapunov coefficient $l_1$ of the three-variable system at its
Hopf points near $\mathrm M$ and $\mathrm N$, computed for $\eps=10^{-3}$ and $10^{-4}$ and the
six values of $\eps'$, has the sign of $A$ in all $24$ cases; as $\eps$
decreases, $\sqrt\eps\,l_1$ converges, as expected at a singular Hopf
point. The zeros of $A$ in $\eps'$, at $0.388375$ for $\mathrm M$ and $0.028855$ for
$\mathrm N$, compare with the sign changes of $l_1$ at $0.388740$ and $0.028849$ for
$\eps=10^{-4}$, and at $0.388411$ and $0.028854$ for $\eps=10^{-5}$; at $\mathrm M$
the gap shrinks tenfold with $\eps$. Simulations for $\eps=10^{-3}$ confirm
the criticality at $\mathrm M$: for $\eps'=1$ the relaxation oscillation persists
down to $f=0.7146225674$, which is $9.4\cdot10^{-5}$ below the Hopf point, and
for $\eps'=10^{-2}$ and $10^{-3}$ small stable cycles grow and jump to the
relaxation oscillation at $1.7\cdot10^{-5}$ and $9.1\cdot10^{-6}$ above it.

\paragraph{Random networks of rank two.}
We drew networks of five to seven reactions on $\ce X,\ce Y,\ce Z$, with
reactant complexes of molecularity at most two, or three if they contain
$2\ce X$, and products of molecularity at most three, always including one
autocatalytic step such as $2\ce X\to3\ce X$; the changes of $\ce X$ and $\ce Y$
are fast and those of $\ce Z$ slow, and we kept networks whose fast part has
rank two. At a random rational point we sampled positive fluxes in
$\{Lv=0,\ gv=0\}$, located $\det M=0$ by bisection between fluxes of opposite
sign of $\det M$, and kept the point if (C1) and (C2) held. The rate
constants $k_r=v_r/x^{y_r}$, rounded to rationals, with the rate constant of
the reaction of largest $|g_rv_r|$ as unfolding parameter, give a canard point
near the chosen one, which was located by exact real root isolation. In
$1\,527$ sampled networks this produced $30$ canard points, $14$
supercritical and $16$ subcritical, with the nonzero eigenvalue of $M$
negative at $28$ of them and positive at $2$. At all $30$ the two
implementations of $A$ agree to $2\cdot10^{-14}$. A Hopf point of the
three-dimensional system near the canard point was found at $58$ of the
$60$ pairs of a canard point and $\eps\in\{10^{-3},10^{-4}\}$, and at each of
them $l_1$ has the sign of $A$. The corrections through $M^\sharp$ change the
sign of $A$ at $8$ of the $30$ canard points, and the transverse term alone at
$5$.

\paragraph{Intrinsic criterion.}
\citet{DeMaesschalckDoanWynen2021} decide the criticality of a planar
slow--fast Hopf point by quantities $V(\mathcal G)(p)$ and $\sigma$ that do
not depend on coordinates. For a planar system of the form (3.1) of
\citet{KrupaSzmolyan2001relax} with Taylor coefficients $F_{zz}$, $F_u$,
\dots, $G_{zz}$ at the fold, a direct computation in the coordinates
$(z,u)$, with $Z=\partial_z$ and the Euclidean metric, gives
\[
  V(\mathcal G)(p)=F_uG_z,\qquad \sigma=\frac{\mathcal Q}{2F_{zz}^2},\qquad
  A=\frac{4\sigma}{\sqrt{-V(\mathcal G)(p)}} .
\]
Indeed, parametrising the critical curve by $u=\varphi(z)$, with
$\varphi''(0)=-F_{zz}/F_u$, their function is $\mathcal G=\hat G\,(F_uF_{zz}-F_zF_{zu})$
on the curve, their vector field $V$ is the derivative with respect to
$F_z$ along the curve, and their skewness is $F_{zzz}/F_{zz}^2$; the identities
follow by differentiating twice at the fold. Applied to the reduced system of
Theorem~\ref{thm:rank}, whose coefficients are~\eqref{eq:formsW}, this shows
that the two criteria agree for every rank. At the five planar benchmark
canard points of the planar paper (the autocatalator, the Brøns--Sturis
model, the Brusselator for two parameter values and a slow-predator model),
each written in three linear charts, two of which make the fast direction
oblique, $A=4\sigma/\sqrt{-V(\mathcal G)(p)}$ holds to all printed digits in
all $15$ cases, with $\sigma$ and $V(\mathcal G)(p)$ independent of the
chart.

\begin{table}[ht]
\centering
\caption{Scripts and logs of the Zenodo record (directory
\texttt{code/julia/models} and \texttt{checks}).}\label{tab:scripts}
\footnotesize
\begin{tabular}{lll}
\toprule
Statement & script & log\\
\midrule
Section~\ref{sm:flux} & \texttt{p2\_rankone.jl} & \texttt{p2\_rankone.txt}\\
Section~\ref{sm:one}, Example~\ref{ex:onemode} & \texttt{p2\_modes.jl} & \texttt{p2\_modes.txt}\\
Section~\ref{sm:lock} & \texttt{p2\_signlock.jl} & \texttt{p2\_signlock.txt}\\
Example~\ref{ex:both}, Figure~\ref{fig:explosion} & \texttt{p2\_global.jl} & \texttt{p2\_global.txt}\\
Example~\ref{ex:both}, failure of (G3) & \texttt{p2\_g3\_failure.py} & \texttt{p2\_g3\_failure.txt}\\
Section~\ref{sm:struct}, two species & \texttt{p2\_g2g3.jl} & \texttt{p2\_g2g3.txt}\\
Section~\ref{sm:struct}, general conditions, weaker (S2) & \texttt{p2\_global\_n3.jl} & \texttt{p2\_global\_n3.txt}\\
Section~\ref{sm:eps}, Example~\ref{ex:both} & \texttt{p2\_bt.jl} & \texttt{p2\_bt.txt}\\
Section~\ref{sm:eps}, three species & \texttt{p2\_bt\_n3.jl} & \texttt{p2\_bt\_n3.txt}\\
Section~\ref{sm:mixing} & \texttt{p2\_mixing.jl} & \texttt{p2\_mixing.txt}\\
Example~\ref{ex:structboth}, Section~\ref{sm:struct}, three species & \texttt{p2\_structural\_both.jl} & \texttt{p2\_structural\_both.txt}\\
Proposition~\ref{prop:oreg2}, exact & \texttt{p2\_oregonator2.py} & \texttt{p2\_oregonator2.txt}\\
Proposition~\ref{prop:oreg2}, $l_1$ & \texttt{p2\_oregonator2\_l1.jl} & \texttt{p2\_oregonator2\_l1.txt}\\
Section~\ref{sm:rank}, flux forms, random networks & \texttt{p2\_highrank.jl} & \texttt{p2\_highrank.txt}\\
Example~\ref{ex:oregonator}, Table~\ref{tab:oregonator} & \texttt{wp3\_coefficient.jl} & \texttt{wp3\_coefficient.txt}\\
Section~\ref{sm:rank}, $l_1$ of the Oregonator & \texttt{criticality.jl} & \texttt{criticality.txt}\\
Section~\ref{sm:rank}, simulations of the Oregonator & \texttt{numerics.jl} & \texttt{numerics.txt}\\
Section~\ref{sm:rank}, intrinsic criterion & \texttt{wp3\_dmdw.jl} & \texttt{wp3\_dmdw.txt}\\
\bottomrule
\end{tabular}
\end{table}

\bibliographystyle{plainnat}
\bibliography{refs}

\end{document}